\documentclass[pdflatex,sn-mathphys-num]{sn-jnl}
\usepackage{graphicx}
\usepackage{amsmath,amssymb,amsfonts}
\usepackage{mathtools}
\usepackage{xcolor}
\usepackage{booktabs}
\usepackage{array}
\usepackage{algorithm}
\usepackage{algorithmic}
\usepackage{etoolbox}
\usepackage{xurl}
\usepackage{cleveref}
\providecommand{\burl}[1]{{\urlstyle{sf}\url{#1}}}
\hypersetup{bookmarksdepth=3}
\makeatletter
\patchcmd{\caption}{\refstepcounter\@captype}
  {\@ifundefined{@float@c@\@captype}{\H@refstepcounter\@captype}{\refstepcounter\@captype}}
  {}{\PackageError{template-adaptation}{Caption hyperlink compatibility patch failed}{Check the float and hyperref package versions.}}
\makeatother

\theoremstyle{thmstyleone}
\newtheorem{theorem}{Theorem}[section]
\newtheorem{lemma}[theorem]{Lemma}
\newtheorem{proposition}[theorem]{Proposition}
\newtheorem{corollary}[theorem]{Corollary}
\theoremstyle{thmstyletwo}
\newtheorem{remark}{Remark}
\newtheorem{example}{Example}
\theoremstyle{thmstylethree}
\newtheorem{assumption}{Assumption}

\crefname{assumption}{Assumption}{Assumptions}
\crefname{remark}{Remark}{Remarks}
\crefname{example}{Example}{Examples}
\crefname{algorithm}{Algorithm}{Algorithms}

\DeclareMathOperator*{\argmin}{arg\,min}
\newcommand{\R}{\mathbb{R}}
\newcommand{\cO}{\mathcal{O}}
\newcommand{\ip}[2]{\langle #1,\, #2 \rangle}
\newcommand{\dom}{\operatorname{dom}}
\newcommand{\rred}[1]{#1}
\newcommand{\intdom}[1]{\text{int}(\dom\,#1)}
\begin{document}
\title{Dual-Inertial Bregman Methods  with Energy-Based Safeguarding}

\author*[1]{\fnm{Kangming} \sur{Chen}}\email{kangming@tmu.ac.jp}

\affil*[1]{\orgdiv{Faculty of Economics and Business Administration},
  \orgname{Tokyo Metropolitan University},
  \orgaddress{
    \postcode{192--0397},
    \state{Tokyo},
    \country{Japan}
  }}

\abstract{We study dual-inertial (DI) Bregman proximal methods for convex optimization via an exact representation of the signed energy drift between dual-coordinate and primal affine interpolation. 
For unrestarted DI, bounded one-sided directional curvature gives a last-iterate objective-gap rate $\cO((1+\log N)/N^2)$ after $N$ iterations, which improves to $\cO(N^{-2})$ under weighted summability. The energy balance also provides a computable acceptance test, which we use in two safeguards. The restart safeguard attains a global $\cO(T^{-1})$ rate after $T$ proximal calls, including rejected trials. Under local quadratic growth and kernel regularity, it restarts only finitely often and eventually converges linearly. The non-restarting safeguard preserves the quadratic weight recurrence via dual substeps. Every fallback terminates without triangle scaling, yielding a $\cO(N^{-2})$ rate in outer iterations and an $\cO(T^{-1})$ rate in calls, which improves to $\cO(T^{-2})$ when the cumulative fallback cost is $\cO(N)$. Experiments on relative-entropy regression and nonnegative least squares illustrate the results.
}



\maketitle

\section{Introduction}\label{sec:intro}

We study acceleration of Bregman proximal-type methods for the convex problem
\begin{equation}\label{eq:intro_problem}
    \min_{x\in\dom\omega} F(x),
\end{equation}
where $\omega$ is a Legendre distance-generating function (Bregman kernel) and $F$ can be extended-valued.

Bregman divergences~\citep{bregman1967} provide a geometry-adapted alternative to Euclidean quadratic regularization and underlie Bregman proximal methods~\citep{censor1992,chen1993,eckstein1993} and mirror descent~\citep{beck2003mirror}.
Relative smoothness permits Bregman descent without a globally Lipschitz gradient~\citep{bauschke2017,lu2018relatively}, but does not by itself guarantee acceleration. Dragomir et al.~\citep{dragomir2021} show that the $\cO(k^{-1})$ worst-case rate is optimal within their Bregman first-order oracle model.

Our approach introduces inertia through interpolation and reconstruction in dual space. For $x_k,z_k\in \operatorname{int}(\dom\omega)$ and $0<\theta_k\le1$, we first form
\[
w_k=\nabla\omega^*\!\left(
(1-\theta_k)\nabla\omega(x_k)
+\theta_k\nabla\omega(z_k)
\right).
\]
After computing $x_{k+1}$ by a Bregman oracle step at $w_k$, we reconstruct $z_{k+1}$ so that
\[
\nabla\omega(x_{k+1})
=(1-\theta_k)\nabla\omega(x_k)
+\theta_k\nabla\omega(z_{k+1}).
\]
We call this the \emph{dual-inertial} (DI) scheme. The unified splitting $F=f+h$ includes the Bregman proximal point method (DI-BPPA) for $f=0$, $h=F$, the Bregman proximal gradient method (DI-BPG) for a general composite splitting, and the Bregman mirror-descent method (DI-BMD) for $h=0$.

\subsection{Background and Related Work}

One approach to Bregman acceleration imposes scaling conditions on the divergence.
Hanzely et al.~\citep{hanzely2021accelerated} obtain $\cO(k^{-\gamma})$ under a triangle-scaling exponent and also study adaptive variants. 
The unified conjugate framework of Gutman and Pe\~na~\citep{gutman2018unified} uses related scaling conditions. 
Liu et al.~\citep{liu2022dual} study dual Bregman proximal gradient methods for relatively strongly convex composite optimization, with accelerated rates determined by the triangle-scaling exponent. For fully implicit proximal steps, Yan and He~\citep{yan2020bregman} develop accelerated Bregman proximal point methods under a quadratic triangle-scaling condition. 
Yang and Toh develop inexact accelerated BPPA under quadrangle scaling~\citep{yang2022bregman} and inertial BPG under restricted relative smoothness and suitable error control~\citep{yang2026inexact}.

Other accelerated methods use compatibility conditions between the objective and the kernel, or control perturbations in a Lyapunov inequality. 
Chen et al.~\citep{chen2026accmd} derive accelerated mirror descent from a variable--operator splitting with a generalized Cauchy--Schwarz condition. Xu and Chen~\citep{xu2026adaptive} propose adaptive accelerated mirror descent (AAMD) based on accumulated Lyapunov perturbation control. For convex objectives, their homotopy variant AAMD-0 gives a $\cO(k_s^{-2})$ energy bound at stage endpoints $k_s$, under dual relative smoothness, mirror compatibility, and bounded Bregman distances. Merchav et al.~\citep{merchav2026adaptive} develop the adaptive integrated mean framework and the Bregman accelerated proximal gradient method (BAPG). Its Shannon entropy specialization, AEPG, achieves a $\cO((1+\log k)/k^2)$ objective-gap rate for its main sequence $\bar x^k$, under a uniform lower bound on its coordinates in the support of an optimal solution.

Primal--dual geometry underlies dual extrapolation~\citep{nesterov2007dual}, linear coupling~\citep{allen2016optimal}, dual-space preconditioning~\citep{maddison2021dual}, and mirror duality~\citep{kim2023mirror}. 
Defazio~\citep{defazio2019curved} interprets accelerated gradient methods through a Bregman proximal construction based on the objective's conjugate.
Dual  inertial interpolation also appears in our companion work on bilevel variational inequalities~\citep{chen2026bregman}. Fast dual proximal gradient methods~\citep{beck2014fast}, by contrast, accelerate a separate dual optimization problem. Azizian et al.~\citep{azizian2024rate} relate local last-iterate rates of first-order Bregman prox-mappings to geometry, regularity, and sharpness.

\subsection{Contributions and Organization}
In this paper, we analyze the signed energy drift caused by the mismatch between dual-coordinate and primal affine interpolation, and use the resulting energy balance to define a common acceptance test for two safeguards: S-DI, which restarts after a rejected trial, and DI-FS, which completes rejected updates via dual substeps without resetting the weights. 
Table~\ref{tab:summary_guarantees} summarizes our convergence guarantees.
\begin{table}[htbp]
\centering
\caption{Summary of convergence guarantees}
\label{tab:summary_guarantees}
\footnotesize
\setlength{\tabcolsep}{3pt}
\begin{tabular*}{\linewidth}{@{\extracolsep{\fill}}
>{\raggedright\arraybackslash}p{0.08\linewidth}
>{\raggedright\arraybackslash}p{0.28\linewidth}
>{\raggedright\arraybackslash}p{0.10\linewidth}
>{\raggedright\arraybackslash}p{0.12\linewidth}
>{\raggedright\arraybackslash}p{0.22\linewidth}
>{\raggedright\arraybackslash}p{0.12\linewidth}@{}}
\toprule
Method & Key conditions & Boundary $x^*$
& Condition type & Guarantee & Ref. \\
\midrule
DI
& $\lambda L\le1$, $\sup_{k\ge2}\eta_k<\infty$;
  additionally $\sum_{k\ge2}\theta_{k-1}\eta_k<\infty$
  for the sharper rate
& Allowed
& Trajectory
& $\cO((1+\log N)/N^2)$;
  $\cO(N^{-2})$ with summability
& Thm.~\ref{prop:directional_budget} \\
\addlinespace[2pt]
DI
& $\lambda L\le1$, DCM
& Allowed
& Trajectory
& $\cO(N^{-2})$, equivalently $\cO(T^{-2})$
& Cor.~\ref{thm:dcm} \\
\midrule
S-DI
& $\lambda L\le1$
& Allowed
& A priori
& $\min\{\cO(T^{-1}),\allowbreak \cO(K_TT^{-2})\}$
& Thm.~\ref{thm:adaptive_bppa_complexity} \\
\addlinespace[2pt]
S-DI
& $\lambda L<1$, local QG,
  $\omega\in C^2$ near $x^*$,
  $\nabla^2\omega(x^*)\succ0$
& Excluded
& Local structure
& Restart bound; eventual linear convergence
& Thm.~\ref{thm:safeguard_finite_restart};
  Cor.~\ref{cor:bppa_eventual_linear} \\
\addlinespace[2pt]
S-DI
& $\lambda L<1$, local weak sharpness on interior solutions
& Excluded
& Local structure
& Finite restarts; finite optimal identification
& Prop.~\ref{prop:safeguard_weak_sharpness} \\
\midrule
DI-FS
& $\lambda L\le1$
& Allowed
& A priori
& 
  $\cO(N^{-2})$;
  $\cO(T^{-1})$ 
& Thm.~\ref{thm:frac_hybrid}(1) \\
\addlinespace[2pt]
DI-FS
& $\lambda L\le1$,
  $C_N\le\bar cN$ for all $N$
& Allowed
& Cost condition
& $\cO(T^{-2})$
& Thm.~\ref{thm:frac_hybrid}(2) \\
\bottomrule
\end{tabular*}

\smallskip
\begin{minipage}{\textwidth}
\footnotesize
\textbf{Notation:} $N$: outer iterations; $T$: total proximal calls, including rejected trials; $C_N$: additional fallback calls through DI-FS outer state $N$; $K_T$: nonempty S-DI epochs after $T$ calls; $\lambda$: stepsize; $L$: relative smoothness constant; $\eta_k$: one-sided directional curvature bound; DCM: directional curvature monotonicity; QG: quadratic growth. \emph{Boundary $x^*$} indicates whether $x^*\in\partial\dom\omega$ is permitted.
\end{minipage}
\end{table}

The main contributions are as follows.
\begin{enumerate}
\item We establish an exact signed connection-drift representation for dual-inertial Bregman methods, yielding last-iterate $\cO((1+\log N)/N^2)$ objective-gap rate after $N$ iterations  under bounded one-sided directional curvature and $\cO(N^{-2})$ under weighted summability.

\item We propose an energy-based restart safeguard (S-DI) that guarantees a global $\cO(T^{-1})$ objective-gap bound after $T$ proximal calls under the standard relative-smoothness stepsize bound. Under a stepsize bound, local quadratic growth at an interior minimizer, and local $C^2$ kernel regularity with positive definite Hessian, we derive an explicit restart bound and prove eventual linear convergence.  Without quadratic growth or $C^2$ kernel regularity, local weak sharpness on interior solutions still ensures finite restarts and finite identification of an optimal point.

\item We develop a non-restarting safeguard (DI-FS) that preserves the quadratic weight recurrence via dual substeps. Each fallback terminates without triangle scaling, even for boundary minimizers. We obtain a global $\cO(N^{-2})$ objective-gap bound in outer iterations and a $\cO(T^{-1})$ bound in total calls. The call bound improves to $\cO(T^{-2})$ if the cumulative fallback cost is $\cO(N)$. 
\end{enumerate}

The rest of this paper is organized as follows. Section~\ref{sec:prelim} reviews the preliminaries, and Section~\ref{sec:algorithms} introduces the unified DI framework. Section~\ref{sec:analysis} establishes the connection-drift energy balance and conditional convergence rates. Section~\ref{sec:adaptive_di} presents the restart-safeguarded scheme and its global and local convergence analysis, while  Section~\ref{ssec:fractional_fallback} develops the non-restarting safeguard based on dual substeps and analyzes its oracle complexity. Section~\ref{sec:experiments_numeric} reports numerical experiments, and Section~\ref{sec:conclusion} concludes the paper.

\section{Preliminaries}
\label{sec:prelim}

\subsection{Bregman Geometry and Basic Identities}

\begin{assumption}\label{ass:omega}
The function $\omega:\R^n\to\R\cup\{+\infty\}$ is proper, closed, Legendre~(see, e.g.,~\cite{rockafellar1970convex}), and supercoercive.
\end{assumption}

In this work, we write \(C:=\operatorname{int}(\dom\omega)\). The Bregman divergence generated by~$\omega$ is defined by the function $D_\omega:\dom\omega\times C\to \mathbb{R}$:
\begin{equation}\label{eq:bregman}
D_\omega(x,y)
:=
\omega(x)-\omega(y)
-\ip{\nabla\omega(y)}{x-y},
\qquad (x,y)\in\dom\omega\times C.
\end{equation}
By convexity, $D_\omega(x,y)\ge0$ and since $\omega$ is Legendre, equality holds if and only if~$x=y$.

The Fenchel conjugate of $\omega$ is
$\omega^*(v):=
\sup_{x\in\R^n}
\{\ip{v}{x}-\omega(x)\}$.
Under Assumption~\ref{ass:omega}, supercoercivity gives $\dom\omega^*=\R^n$. By Legendre duality, $\omega^*$ is also Legendre and the gradient maps are mutually inverse (see \citep[Fact~2.9]{bauschke1997legendre} and~\citep[Theorem~26.5]{rockafellar1970convex}):
\[
\nabla\omega:C\to\R^n,\qquad
\nabla\omega^*:\R^n\to C,\qquad
(\nabla\omega)^{-1}=\nabla\omega^*.
\]

The corresponding divergences satisfy
\begin{equation}\label{eq:bregman_dual}
D_\omega(x,y)
=
D_{\omega^*}\bigl(\nabla\omega(y),\nabla\omega(x)\bigr),
\qquad x,y\in C.
\end{equation}
Differentiation in the first argument gives
\[
\nabla_x D_\omega(x,y)
=
\nabla\omega(x)-\nabla\omega(y),\qquad x,y\in C.
\]

We recall the three-point identity for Bregman divergences~\citep{bauschke1997legendre}.
\begin{lemma}[Three-Point Identity]\label{lem:three_point}
For any $x\in\dom\omega$ and $y,z\in\intdom{\omega}$:
\begin{equation}\label{eq:three_point}
  D_\omega(x,z) \;=\; D_\omega(x,y) + D_\omega(y,z) +
  \ip{\nabla\omega(y)-\nabla\omega(z)}{x-y}.
\end{equation}
\end{lemma}

\begin{lemma}[Bregman sublevel sets]
\label{lem:bregman_sublevel_compact}
Under Assumption~\ref{ass:omega}, for every $q\in C$ and $r\ge0$,
the set $\{z\in C:D_\omega(q,z)\le r\}$ is compact, and its
image under $\nabla\omega$ is bounded.
\end{lemma}
\begin{proof}
Choose $\delta>0$ such that $\overline B_\delta(q)\subset C$.
For $s=\nabla\omega(z)$, Fenchel's inequality gives
\[
D_\omega(q,z)
=\omega^*(s)-\langle q,s\rangle+\omega(q)
\ge\delta\|s\|-\max_{\overline B_\delta(q)}\omega+\omega(q).
\]
Hence $\nabla\omega(z)$ is uniformly bounded whenever
$D_\omega(q,z)\le r$.
By Legendre duality, this image equals
\[
\{s\in\R^n:\omega^*(s)-\langle q,s\rangle+\omega(q)\le r\},
\]
which is closed by continuity of $\omega^*$, so compact.
Its image under the continuous map $\nabla\omega^*$ is the
original sublevel set, and is therefore compact.
\end{proof}

\subsection{\rred{Objective Assumptions and Proximal Well-Posedness}}

We first introduce the assumptions on the objective function, followed by the composite structure required for the proximal gradient method.

\begin{assumption} \label{ass:objective}
$F$ is proper, closed, and convex, with $\dom F\cap C\ne\varnothing$. Moreover, problem~\eqref{eq:intro_problem} has a minimizer $x^*\in\dom F\cap\dom\omega$.
\end{assumption}

\begin{assumption}[Composite structure and relative smoothness]
\label{ass:smad}
For the Bregman proximal gradient method, suppose that
$F=f+h$, where $h$ is proper, closed, and convex, and
$f$ is finite, convex, and differentiable on $C$.
Furthermore, $f$ is $L$-smooth relative to $\omega$
for some $L\ge 0$, namely,
\begin{equation}\label{eq:rel_smooth_def}
f(x)\le f(y)+\ip{\nabla f(y)}{x-y}+L D_\omega(x,y),
\qquad x,y\in C.
\end{equation}
Equivalently, $L\omega-f$ is convex on $C$.
\end{assumption}
Throughout the remainder of the paper, Assumptions~\ref{ass:omega}, \ref{ass:objective}, and~\ref{ass:smad} are our standing assumptions.
Strong convexity is taken with respect to the Euclidean norm. A function $F$ is $m_F$-strongly convex, with $m_F>0$, if $F-\frac{m_F}{2}\|\cdot\|^2$ is convex. In this case,
\[
F(x)\ge F(y)+\langle g,x-y\rangle
+\frac{m_F}{2}\|x-y\|^2,
\qquad x,y\in\dom F,\quad g\in\partial F(y).
\]
For relative strong convexity with respect to $\omega$, the quadratic term is replaced by $m_{\mathrm{rel}}D_\omega(x,y)$, where $m_{\mathrm{rel}}>0$ and $y\in C$.

The following lemma ensures that the Bregman proximal subproblems used below have unique solutions in $C$.
\begin{lemma}\label{lem:prox_interior}
Under Assumption~\ref{ass:omega}, let
$\varphi:\R^n\to\R\cup\{+\infty\}$ be proper, closed, and convex,
with $\dom\varphi\cap C\ne\varnothing$.
For any $w\in C$ and $\lambda>0$, the problem
\[
\min_x\left\{\varphi(x)+\lambda^{-1}D_\omega(x,w)\right\}
\]
has a unique solution $x^+\in\dom\varphi\cap C$.
Moreover, $ \lambda^{-1} \bigl(\nabla\omega(w)-\nabla\omega(x^+)\bigr) \in\partial\varphi(x^+)$.
For fixed $\varphi$ and $\lambda$, the unique minimizer of
$
\omega(x)+\lambda\varphi(x)-\langle t,x\rangle
$
depends continuously on $t\in\R^n$.
\end{lemma}

\begin{proof}
Existence, uniqueness, and the optimality condition follow
from~\citep[Propositions~3.21(vii), 3.22(ii)(a),(d),
and 3.23(v)(b)]{bauschke2003bregman}.

For continuity, set $\psi:=\omega+\lambda\varphi$,
which is proper, closed, and convex.
For each $t\in\R^n$, the first part with
$w=\nabla\omega^*(t)$ gives a unique minimizer $p(t)$ of
$\psi-\langle t,\cdot\rangle$.
Thus $\psi^*$ is finite on $\R^n$, and from Fenchel conjugacy, we have
\[
\partial\psi^*(t)
=\argmin_x\{\psi(x)-\langle t,x\rangle\}
=\{p(t)\}.
\]
Hence $\psi^*$ is continuously differentiable and
$p(t)=\nabla\psi^*(t)$ is continuous
\citep[Theorem~25.1 and Corollary~25.5.1]
{rockafellar1970convex}.
\end{proof}
\section{Dual-Inertial Bregman Methods}\label{sec:algorithms}

\subsection{Dual-Inertial Interpolation}\label{ssec:di_principle}
For $x,z\in C$ and $\tau\in\mathbb{R}$, define
\begin{equation}\label{eq:dual_coordinate_map}
T_\tau(x,z)
:=
\nabla\omega^*
\bigl((1-\tau)\nabla\omega(x)+\tau\nabla\omega(z)\bigr).
\end{equation}
Supercoercivity in Assumption~\ref{ass:omega} gives
$\dom\omega^*=\R^n$. Thus $T_\tau(x,z)$ is well-defined in $C$
for every $\tau\in\R$.
Moreover,
\begin{equation}\label{eq:dual_coordinate_pullback}
\nabla\omega(T_\tau(x,z))
=
(1-\tau)\nabla\omega(x)+\tau\nabla\omega(z).
\end{equation}
So $T_\tau$ performs affine interpolation in dual coordinates when $\tau\in[0,1]$, and extrapolation otherwise.

\begin{lemma}\label{lem:dual_interpolation_composition}
For $x,z\in C$ and $\alpha,\beta,\gamma\in\mathbb R$,
\begin{equation}\label{eq:dual_interpolation_composition}
T_\gamma(T_\alpha(x,z),T_\beta(x,z))
=T_{(1-\gamma)\alpha+\gamma\beta}(x,z).
\end{equation}
In particular,
\begin{align}
T_\beta(x,T_\alpha(x,z))&=T_{\alpha\beta}(x,z),\qquad
T_\beta(T_\alpha(x,z),z)=T_{\alpha+\beta(1-\alpha)}(x,z),
\label{eq:dual_interpolation_composition_special}\\
T_{1/\theta}(x,T_\theta(x,z))&=z,\qquad \theta\ne0.
\label{eq:dual_interpolation_inverse}
\end{align}
\end{lemma}
\begin{proof}
Set $u=\nabla\omega(x)$ and $s=\nabla\omega(z)$. Applying \eqref{eq:dual_coordinate_pullback} twice gives
\[
\begin{aligned}
\nabla\omega(T_\gamma(T_\alpha(x,z),T_\beta(x,z)))
&=(1-\gamma)((1-\alpha)u+\alpha s)
  +\gamma((1-\beta)u+\beta s)\\
&=(1-((1-\gamma)\alpha+\gamma\beta))u+ ((1-\gamma)\alpha+\gamma\beta) s,
\end{aligned}
\]
which proves \eqref{eq:dual_interpolation_composition}. Since $T_0(x,z)=x$ and $T_1(x,z)=z$, taking the parameter triples $(0,\alpha,\beta)$ and $(\alpha,1,\beta)$ gives the two identities in~\eqref{eq:dual_interpolation_composition_special}. Taking $\alpha=\theta$ and $\beta=1/\theta$ in the first one gives $T_1(x,z)=z$, proving~\eqref{eq:dual_interpolation_inverse}.
\end{proof}

For $x_k,z_k\in C$ and $0<\theta_k\le1$, define the interpolation center by
\begin{equation}\label{eq:vk_param}
w_k=T_{\theta_k}(x_k,z_k).
\end{equation}
After computing $x_{k+1}$, choose the auxiliary point $z_{k+1}$ so that
\begin{equation}\label{eq:coupling_reconstruction}
x_{k+1}=T_{\theta_k}(x_k,z_{k+1}).
\end{equation}
By \eqref{eq:dual_interpolation_inverse}, this uniquely determines
\begin{equation}\label{eq:auxiliary_reconstruction}
z_{k+1}=T_{1/\theta_k}(x_k,x_{k+1}).
\end{equation}
Thus the auxiliary point is obtained by dual coordinate extrapolation from $x_k$ through $x_{k+1}$. For the Euclidean kernel $\omega(x)=\frac12\|x\|^2$, it reduces to
$z_{k+1}
=
x_k+\frac{1}{\theta_k}(x_{k+1}-x_k)$.

For the convergence analysis, we choose the interpolation parameters through the accelerated weight recurrence~\citep{beck2009fista,monteiro2013accelerated}:
\begin{equation}\label{eq:param}
A_k=(1-\theta_k)A_{k+1},
\qquad
A_{k+1}\theta_k^2=\lambda,
\qquad
A_1>0,\quad \lambda>0.
\end{equation}
Equivalently, setting $t_k:=\sqrt{A_k/\lambda}$ gives
\[
t_{k+1}=\frac{1+\sqrt{1+4t_k^2}}{2},
\qquad
\theta_k=\frac{1}{t_{k+1}},
\qquad
t_1=\sqrt{A_1/\lambda}.
\]
Since $1/2\le t_{k+1}-t_k\le1$, we have $A_k=\Theta(k^2)$, $\theta_k=\Theta(k^{-1})$, and $\theta_k\le\theta_{k-1}$. Moreover, $t_k\ge t_1+(k-1)/2$ gives, for $N\ge2$,
\[
\sum_{k=2}^N\theta_{k-1}
\le2\sum_{j=1}^{N-1}\frac1j
\le2(1+\log N).
\]
We also define the weight increment
\begin{equation}\label{eq:weight_increment}
a_{k+1}
:=
A_{k+1}-A_k
=
A_{k+1}\theta_k
=
\frac{\lambda}{\theta_k}.
\end{equation}
With the standard initialization $A_1=\lambda$, equivalently $t_1=1$, we further have
\begin{equation}\label{eq:weight_lower_bound}
A_k \ge \frac{\lambda(k+1)^2}{4}.
\end{equation}

\subsection{The DI Algorithm and Its Special Cases}

For $w\in C$ and $\lambda>0$, define the Bregman proximal update
\begin{equation}\label{eq:prox_map}
P_\lambda(w)
:=\argmin_{y\in\dom\omega}
\left\{
\langle\nabla f(w),y\rangle+h(y)
+\lambda^{-1}D_\omega(y,w)
\right\}.
\end{equation}
Algorithm~\ref{alg:di_framework} uses this map at the
dual-interpolated center $w_k$.

\begin{algorithm}[H]
\caption{Dual-inertial Bregman algorithm (DI)}
\label{alg:di_framework}
\begin{algorithmic}[1]
\REQUIRE $x_1=z_1\in\dom F\cap C$;
$\lambda>0$; $0<\theta_k\le1$.
\FOR{$k=1,2,\ldots$}
\STATE $w_k=T_{\theta_k}(x_k,z_k)$
\STATE $x_{k+1}=P_\lambda(w_k)$
\STATE $z_{k+1}=T_{1/\theta_k}(x_k,x_{k+1})$
\ENDFOR
\end{algorithmic}
\end{algorithm}

Taking $f=0$ and $h=F$ gives exact DI-BPPA with $L=0$.
The general composite case gives DI-BPG, while $h=0$ gives
DI-BMD with
\[
x_{k+1}
=\nabla\omega^*\bigl(
\nabla\omega(w_k)-\lambda\nabla f(w_k)
\bigr).
\]
The rate analysis uses \eqref{eq:param} and $\lambda L\le1$, and the latter is automatic for DI-BPPA.

\begin{corollary}[Well-posedness]\label{cor:di_wellposed}
Under Assumptions~\ref{ass:omega}, \ref{ass:objective} and~\ref{ass:smad}, every finite iteration of Algorithm~\ref{alg:di_framework} is well-defined. Each subproblem has a unique solution, and $x_k\in\dom F\cap C$, $w_k,z_k\in C$.
\end{corollary}
\begin{proof}
Well-posedness follows from Lemma~\ref{lem:prox_interior}, applied to $h+\langle\nabla f(w_k),\cdot\rangle$, and $T_\tau(C\times C)\subset C$ for every $\tau\in\R$.
\end{proof}

\subsection{Common Oracle Relation and Dual Update}
\begin{proposition}
\label{prop:common_dual_update}
Consider Algorithm~\ref{alg:di_framework} with $\lambda>0$ and $0<\theta_k\le1$, and define
\begin{equation}\label{eq:oracle_residual}
g_{k+1}
:=
\lambda^{-1}
\bigl(
\nabla\omega(w_k)-\nabla\omega(x_{k+1})
\bigr).
\end{equation}
Then
\[
g_{k+1}\in\nabla f(w_k)+\partial h(x_{k+1}).
\]
Moreover,
\begin{equation}\label{eq:auxiliary_dual_update}
\nabla\omega(z_{k+1})-\nabla\omega(z_k)
=
-\frac{\lambda}{\theta_k}g_{k+1}.
\end{equation}
\end{proposition}
\begin{proof}
The oracle inclusion follows from the first-order optimality condition of the subproblem. 
Subtracting the dual-coordinate identities
\eqref{eq:vk_param} and \eqref{eq:coupling_reconstruction}, and using \eqref{eq:oracle_residual}, proves~\eqref{eq:auxiliary_dual_update}.
\end{proof}
\begin{remark}
For $k\ge2$, the interpolation center can be written as
\[
w_k=T_{-\alpha_k}(x_k,x_{k-1}),
\qquad
\alpha_k=
\frac{\theta_k(1-\theta_{k-1})}{\theta_{k-1}}.
\]
Indeed, by Lemma~\ref{lem:dual_interpolation_composition},
\[
w_k
=T_{\theta_k}\!\left(
x_k,T_{1/\theta_{k-1}}(x_{k-1},x_k)
\right)
=T_{-\alpha_k}(x_k,x_{k-1}).
\]
With $x_0=x_1$ and $\alpha_1=0$, the same formula holds at initialization. So DI methods can be implemented using the two most recent iterates, without explicitly storing $z_k$.
\end{remark}

\section{Connection Drift and Convergence Rates}
\label{sec:analysis}
We derive an exact energy balance for the DI scheme and identify the signed connection drift. Bounds on its cumulative contribution then yield last-iterate convergence rates.
Unless stated otherwise, we consider exact DI updates
with fixed $\lambda>0$, weights satisfying~\eqref{eq:param},
and $x_1=z_1\in\dom F\cap C$.

\subsection{Energy Balance and Connection Drift}
\label{ssec:connection_drift}

For the objective $F$, define the Lyapunov function
\begin{equation}\label{eq:lyap_bppa}
  E_k \;=\; A_k\bigl(F(x_k)-F(x^*)\bigr)
           \;+\; D_\omega(x^*, z_k).
\end{equation}
By nonnegativity of the Bregman divergence, $E_k\ge 0$ for all $k$.

Define the signed interpolation defect
\begin{equation}\label{eq:jensen_gap}
\Delta_k:=\ip{g_{k+1}}
{(1-\theta_k)x_k+\theta_kz_{k+1}-x_{k+1}},
\end{equation}
which measures the mismatch between primal affine and dual-coordinate interpolation in the direction $g_{k+1}$.

\begin{theorem}[Exact energy balance]
\label{thm:di_bppa}\label{rem:connection_slack}
Let $D_f(y,w)=f(y)-f(w)-\langle\nabla f(w),y-w\rangle$ and define
\[
\begin{aligned}
U_k&:=A_k\bigl[F(x_k)-F(x_{k+1})
-\langle g_{k+1},x_k-x_{k+1}\rangle+D_f(x_{k+1},w_k)\bigr],\\
V_k&:=F(x^*)-F(x_{k+1})
-\langle g_{k+1},x^*-x_{k+1}\rangle+D_f(x_{k+1},w_k),\\
S_k&:=U_k+(A_{k+1}-A_k)V_k.
\end{aligned}
\]
Then $U_k,V_k,S_k\ge0$ and, for every $k\ge1$,
\begin{equation}\label{eq:exact_bppa_energy}
E_{k+1}-E_k
=-D_\omega(z_{k+1},z_k)-A_{k+1}\Delta_k-S_k
+A_{k+1}D_f(x_{k+1},w_k).
\end{equation}
In particular,
\begin{equation}\label{eq:lyap_descent_bppa}
E_{k+1}-E_k
\le-D_\omega(z_{k+1},z_k)-A_{k+1}\Delta_k
+A_{k+1}L D_\omega(x_{k+1},w_k).
\end{equation}
\end{theorem}
\begin{proof}
The oracle relation gives $\xi_{k+1}:=g_{k+1}-\nabla f(w_k)\in\partial h(x_{k+1})$. Since $F$ and $h$ are closed and convex, $f=F-h$ is radially continuous at $x^*$ along a segment to a point in $\dom F\cap C$. Thus the convex support inequality for $f$ extends to $x^*$. For $y\in\{x_k,x^*\}$, convexity gives
\[
\begin{aligned}
&F(y)-F(x_{k+1})-\langle g_{k+1},y-x_{k+1}\rangle
+D_f(x_{k+1},w_k)\\
&=D_f(y,w_k)+h(y)-h(x_{k+1})
-\langle\xi_{k+1},y-x_{k+1}\rangle\ge0.
\end{aligned}
\]
Hence $U_k,V_k,S_k\ge0$. 
Since $A_k=(1-\theta_k)A_{k+1}$, the definitions of the remainders can be rearranged as
\begin{equation}\label{eq:common_objective_telescope}
\begin{aligned}
&A_{k+1}(F(x_{k+1})-F(x^*))-A_k(F(x_k)-F(x^*))\\
&=A_{k+1}\ip{g_{k+1}}
{x_{k+1}-(1-\theta_k)x_k-\theta_kx^*}
+A_{k+1}D_f(x_{k+1},w_k)-S_k.
\end{aligned}
\end{equation}
The three-point identity \eqref{eq:three_point}, the auxiliary update \eqref{eq:auxiliary_dual_update}, and $A_{k+1}\theta_k^2=\lambda$ give
\begin{equation}\label{eq:common_bregman}
D_\omega(x^*,z_{k+1})-D_\omega(x^*,z_k)
=-D_\omega(z_{k+1},z_k)
+A_{k+1}\theta_k\langle g_{k+1},x^*-z_{k+1}\rangle.
\end{equation}
Adding \eqref{eq:common_bregman} to \eqref{eq:common_objective_telescope} cancels the $x^*$ terms and proves \eqref{eq:exact_bppa_energy}. Finally, $S_k\ge0$ and $D_f(x_{k+1},w_k)\le L D_\omega(x_{k+1},w_k)$ imply \eqref{eq:lyap_descent_bppa}.
\end{proof}

To analyze the signed term $\Delta_k$, set
\[
H:=\omega^*,\qquad
u_k:=\nabla\omega(x_k),\qquad
s_k:=\nabla\omega(z_k).
\]

For $x,z\in C$ and $0<\theta\le1$, define the connection correction
\begin{equation}\label{eq:connection_R}
\mathcal R_\theta(x,z)
:=D_\omega(x,z)-\theta^{-2}D_\omega(x,T_\theta(x,z)).
\end{equation}

For $u=\nabla\omega(x)$ and $s=\nabla\omega(z)$, Bregman conjugacy gives
$
\mathcal R_\theta(x,z)
=D_H(s,u)-\theta^{-2}D_H\bigl(u+\theta(s-u),u\bigr).
$
Note that for a quadratic kernel, dual and primal affine interpolation coincide, so $\Delta_k=0$ and $\mathcal R_\theta\equiv0$.

Define 
\[
r_k:=\mathcal R_{\theta_k}(x_k,z_k),
\qquad
\Phi_k:=E_k-r_k,
\]
and
\begin{equation}\label{eq:connection_gamma}
\Gamma_{k+1}
:=\mathcal R_{\theta_k}(x_k,z_{k+1})
-\mathcal R_{\theta_{k+1}}(x_{k+1},z_{k+1}).
\end{equation}

Set
\begin{equation}\label{eq:connection_dissipation}
\widehat Q_k
:=\frac{D_\omega(x_{k+1},w_k)-\lambda D_f(x_{k+1},w_k)}
{\theta_k^2},
\qquad
Q_k:=\frac{1-\lambda L}{\theta_k^2}
D_\omega(x_{k+1},w_k).
\end{equation}
Relative smoothness gives $\widehat Q_k\ge Q_k$,
with both nonnegative when $\lambda L\le1$.
\begin{proposition}
\label{prop:connection_edge}
For every $k\ge1$, we have
\begin{equation}\label{eq:exact_connection_recursion}
\Phi_{k+1}+\widehat Q_k+S_k=\Phi_k+\Gamma_{k+1}.
\end{equation}
Consequently,
\begin{equation}\label{eq:connection_recursion}
\Phi_{k+1}+Q_k
\le\Phi_k+\Gamma_{k+1}.
\end{equation}
For DI-BPPA, $\widehat Q_k=Q_k$ with $L=0$.
\end{proposition}
\begin{proof}
The DI updates give
$w_k=T_{\theta_k}(x_k,z_k)$ and
$x_{k+1}=T_{\theta_k}(x_k,z_{k+1})$.
Hence
\begin{equation}\label{eq:connection_balance}
\begin{aligned}
&r_k-\mathcal R_{\theta_k}(x_k,z_{k+1})\\
&=D_\omega(x_k,z_k)-D_\omega(x_k,z_{k+1})
-\theta_k^{-2}
\bigl[D_\omega(x_k,w_k)-D_\omega(x_k,x_{k+1})\bigr]\\
&=D_\omega(z_{k+1},z_k)
-\theta_k^{-2}D_\omega(x_{k+1},w_k)
+A_{k+1}\Delta_k.
\end{aligned}
\end{equation}
The first equality follows from the definition of
$\mathcal R_\theta$.
The second uses the three-point identity, the dual updates
in Proposition~\ref{prop:common_dual_update},
$A_{k+1}\theta_k^2=\lambda$, and
the definition of $\Delta_k$.
Substitution into \eqref{eq:exact_bppa_energy} proves
\eqref{eq:exact_connection_recursion}.
Finally, $S_k\ge0$ and $\widehat Q_k\ge Q_k$ imply
\eqref{eq:connection_recursion}.
\end{proof}

Summing the corrected energy inequality gives the following objective-gap bound in terms of the cumulative drift and terminal correction.

\begin{corollary}
\label{cor:drift_budget}
Assume $\lambda L\le1$ and define
\[
\mathcal B_N:=r_N+\sum_{k=2}^N\Gamma_k,
\]
with $Q_k$ as in \eqref{eq:connection_dissipation}. Then, for every $N\ge2$,
\begin{equation}\label{eq:drift_budget_finite}
0\le A_N\bigl(F(x_N)-F(x^*)\bigr)
\le E_N+\sum_{k=1}^{N-1}Q_k
\le E_1+\mathcal B_N.
\end{equation}
\end{corollary}

\begin{proof}
Since $x_1=z_1$, we have $r_1=0$ and $\Phi_1=E_1$. Summing \eqref{eq:connection_recursion} and substituting $\Phi_N=E_N-r_N$ proves \eqref{eq:drift_budget_finite}. The objective bound follows from $Q_k\ge0$ and $E_N\ge A_N(F(x_N)-F(x^*))$.
In particular,  if $\mathcal B_N\le b_N$ and $b_N/A_N\to0$, then $F(x_N)\to F(x^*)$.
\end{proof}

\subsection{Rates under One-Sided Curvature Control}
\label{ssec:nearlog}
We express the drift through the directional curvature of $H$ along the generated dual segments. The representation below uses a nonnegative Peano kernel~\citep[Chapter~22]{powell1981approximation}.

\begin{lemma}[Exact drift representation]
\label{lem:connection_integral}
For $k\ge2$, set $d_k:=s_k-u_{k-1}$. Fix such a $k$ and write
$p=\theta_{k-1}$, $u=u_{k-1}$, and $d=d_k$. For $q\in(0,1]$,
define $\rho=p+q(1-p)$ and
\[
\Gamma_{k,q}:=\mathcal R_p(x_{k-1},z_k)-\mathcal R_q(x_k,z_k).
\]
If $H$ is $C^3$ on a neighborhood of $[u,u+\rho d]$, then
\begin{equation}\label{eq:gamma_integral}
\Gamma_{k,q}=\int_0^\rho W_{p,q}(t)\nabla^3H(u+td)[d,d,d]\,dt,
\end{equation}
where the nonnegative weight $W_{p,q}$ is given in \eqref{eq:gamma_weight} and satisfies
\begin{equation}\label{eq:gamma_mass}
I(p,q):=\int_0^\rho W_{p,q}(t)\,dt
=\frac{1-p}{6}\bigl[p(2-p)+q(1-p)^2\bigr]
\le\frac{2p+q}{6}.
\end{equation}
In particular, $\Gamma_{k,\theta_k}=\Gamma_k$ and $\Gamma_{k,1}=\Gamma_k+r_k$.
\end{lemma}
\begin{proof}
  The weight and the proof are given in Appendix~\ref{app:drift_representation}.
\end{proof}

\begin{theorem}[One-sided curvature control]\label{prop:directional_budget}
Assume $\lambda L\le1$. For every $k\ge2$, suppose $H$ is $C^3$ on an open neighborhood of the full segment $[u_{k-1},s_k]$, and define
\[
\eta_k:=\sup_{0\le t\le1}
\bigl[\nabla^3H(u_{k-1}+td_k)[d_k,d_k,d_k]\bigr]_+,
\]
where $[a]_+:=\max\{a,0\}$ denotes the positive part of $a$. Then, for $N\ge2$,
\begin{equation}\label{eq:one_sided_directional_budget}
E_N+\sum_{k=1}^{N-1}Q_k
\le E_1+\sum_{k=2}^{N-1}I(\theta_{k-1},\theta_k)\eta_k
+\frac{1-\theta_{N-1}}6\eta_N.
\end{equation}
If $\sup_{k\ge2}\eta_k<\infty$, the last-iterate objective gap is $\cO((1+\log N)/N^2)$. If, in addition, $\sum_{k=2}^\infty\theta_{k-1}\eta_k<\infty$, the rate improves to $\cO(N^{-2})$.
\end{theorem}
\begin{proof}
Since $W_{p,q}\ge0$, Lemma~\ref{lem:connection_integral} gives, for $2\le k<N$,
\[
\Gamma_k\le I(\theta_{k-1},\theta_k)\eta_k
\le(2\theta_{k-1}+\theta_k)\eta_k/6.
\]
For the final terms, let $p=\theta_{N-1}$. Since $\mathcal R_1=0$, the same lemma with auxiliary parameter $q=1$ yields
\[
\Gamma_N+r_N
=\mathcal R_p(x_{N-1},z_N)-\mathcal R_1(x_N,z_N)\le I(p,1)\eta_N=\frac{1-p}{6}\eta_N.
\]
Substitution into \eqref{eq:drift_budget_finite} proves \eqref{eq:one_sided_directional_budget}. The rates follow from $A_N=\Theta(N^2)$, $\theta_k=\Theta(k^{-1})$, and $\theta_k\le\theta_{k-1}$ in Section~\ref{ssec:di_principle}.
\end{proof}

Appendix~\ref{app:drift} gives a $C^2$ alternative based on directional curvature monotonicity, together with sufficient cone conditions for it.

In the following example, DI-BPPA can be rewritten as an accelerated proximal point method in a fixed metric, so no directional curvature assumption is needed.
\begin{example}[Quadratic objectives under Legendre geometry]
\label{rem:kl_quadratic}
Under Assumption~\ref{ass:omega}, consider
\[
\min_{x\in\dom\omega}F(x)
:=\tfrac12x^\top Bx-\langle b,x\rangle,
\qquad 0\ne B=B^\top\succeq0,
\]
and suppose that $\{x\in C:Bx=b\}\ne\varnothing$. Consider exact DI-BPPA with fixed $\lambda>0$ and $x_1=z_1\in C$. The Bregman projection
$
x^\dagger:=\arg\min_{\substack{x\in\dom\omega\\ Bx=b}}
D_\omega(x,x_1)
$
exists uniquely and belongs to $C$. Let $u_1=\nabla\omega(x_1)$ and $u^\dagger=\nabla\omega(x^\dagger)$, and set
$
C_{\rm init}:=A_1D_\omega(x^\dagger,x_1)
+\tfrac12\|u_1-u^\dagger\|_{B^\dagger}^2,
$
where $B^\dagger$ is the Moore--Penrose inverse and $\|v\|_{B^\dagger}^2=v^\top B^\dagger v$. Then, for $N\ge2$,
\[
D_\omega(x^\dagger,x_N)\le C_{\rm init}/A_N,
\quad
F(x_N)-F^*\le C_{\rm init}/(\lambda A_N),
\]
where $F^*=\min_{x\in\dom\omega}F(x)$. Appendix~\ref{app:kl_quadratic} gives the proof.
\end{example}

\section{Restart-Based Safeguard Algorithm}
\label{sec:adaptive_di}
Whereas the results in Section~\ref{sec:analysis} concern the unrestarted  DI trajectory, in this section we use the  energy balance to construct a computable acceptance test that safeguards each DI iteration. 

\subsection{Acceptance Test and Restart Rule}
Consider one trial of Algorithm~\ref{alg:di_framework} with the coupling \eqref{eq:param}. Suppressing the iteration index, write $(A,x,z)$ and $(A^+,x^+,z^+)$ for the current and trial states, where $A>0$, $x\in\dom F\cap C$, and $z\in C$, and set $a=A^+-A$ as in \eqref{eq:weight_increment}. Let $w$ be the trial center and $g$ be defined by \eqref{eq:oracle_residual}. 
By Proposition~\ref{prop:common_dual_update},
$\nabla\omega(z^+)-\nabla\omega(z)=-ag$.
We use the linearization error $D_f$ and $U,V$ from {Theorem~\ref{thm:di_bppa}}. The acceptance margin is
\begin{equation}\label{eq:adaptive_margin}
\mathcal M
:=D_\omega(z^+,z)+A\bigl[F(x)-F(x^+)\bigr]
+a\bigl[\langle g,z^+-x^+\rangle-D_f(x^+,w)\bigr].
\end{equation}
The trial is accepted when
\begin{equation}\label{eq:adaptive_accept_equivalent}
\mathcal M\ge0.
\end{equation}
The test requires neither $x^*$ nor $F(x^*)$ and does not evaluate $F$ at the extrapolated auxiliary points $z$ or $z^+$.
Algorithm~\ref{alg:adaptive_bppa} restarts from the last accepted auxiliary point after a  rejected trial. Each epoch begins with a non-inertial Bregman step from $c\in C$, possibly outside $\dom h$, and its output $p$ initializes the DI state as $(x,z,A)=(p,p,\lambda)$. Here  $T$ counts all Bregman proximal oracle calls, including those producing rejected updates.

\begin{algorithm}[htbp]
\caption{Restart safeguarded dual-inertial Bregman algorithm (S-DI)}
\label{alg:adaptive_bppa}
\begin{algorithmic}[1]
\REQUIRE $c_1\in C$, $\lambda>0$ with
$\lambda L\le1$, and an integer $T\ge1$.
\STATE Set $c=c_1$ and $A=0$; leave $\widehat x$ unset.
\FOR{$n=1,\ldots,T$}
    \IF{$A=0$}
        \STATE $p=P_\lambda(c)$
        \STATE Set $(x,z,A)=(p,p,\lambda)$ and accept $p$.
    \ELSE
        \STATE $a=(\lambda+\sqrt{\lambda^2+4\lambda A})/2$,
        $A^+=A+a$, $\theta=a/A^+$.
        \STATE Generate $(w,x^+,z^+)$ using the update steps of
        Algorithm~\ref{alg:di_framework} at $(x,z)$ with
        step size $\lambda$ and parameter $\theta$.
        \STATE Calculate $g$ by \eqref{eq:oracle_residual} and
        $\mathcal M$ by \eqref{eq:adaptive_margin}.
        \IF{$\mathcal M\ge0$}
            \STATE Set $(x,z,A)=(x^+,z^+,A^+)$
            and accept $x^+$.
        \ELSE
            \STATE Restart: set $(c,A)=(z,0)$.
        \ENDIF
    \ENDIF
    \STATE On acceptance, set $\widehat x=x$ if unset or
    if $F(x)<F(\widehat x)$.
\ENDFOR
\RETURN $\widehat x_T:=\widehat x$.\end{algorithmic}
\end{algorithm}

\begin{proposition}
\label{prop:adaptive_safeguard}
For the preceding trial with any $\lambda>0$, let $E$ and $E^+$ denote the energies \eqref{eq:lyap_bppa} at $(A,x,z)$ and $(A^+,x^+,z^+)$. With $V\ge0$ at $x^*$ from Theorem~\ref{thm:di_bppa}, we have
\begin{equation}\label{eq:adaptive_edge_energy}
E^+-E=-\mathcal M-aV.
\end{equation}
Consequently, \eqref{eq:adaptive_accept_equivalent} is a sufficient condition for $E^+\le E$.
\end{proposition}

\begin{proof}
Using $A=(1-\theta)A^+$ and $a=\theta A^+$, the definitions of $U$, $\Delta$, and $\mathcal M$ give
\[
\mathcal M
=D_\omega(z^+,z)+A^+\Delta+U-A^+D_f(x^+,w).
\]
Substituting this identity and $S=U+aV$ into \eqref{eq:exact_bppa_energy} proves \eqref{eq:adaptive_edge_energy}. Since $a>0$ and $V\ge0$ by Theorem~\ref{thm:di_bppa}, $\mathcal M\ge0$ implies $E^+\le E$.
\end{proof}

\begin{remark}[Geometric interpretation of the margin]
For a DI trial at iteration $k$, let $\mathcal M_k$ denote
the margin in \eqref{eq:adaptive_margin} and set
$
J_k:=\mathcal R_{\theta_k}(x_k,z_{k+1})
     -\mathcal R_{\theta_k}(x_k,z_k).
$
By \eqref{eq:connection_balance},
\begin{equation}\label{eq:margin_decomposition}
\mathcal M_k=\widehat Q_k+U_k-J_k.
\end{equation}
Thus the trial is accepted if and only if
$J_k\le\widehat Q_k+U_k$. By definition,
$J_k=\Gamma_{k+1}+r_{k+1}-r_k$.
\end{remark}
\subsection{Global Oracle Call Bounds}
To combine the energy bounds across epochs, we first estimate
the energy at epoch initialization.
\begin{lemma} 
\label{lem:adaptive_epoch_transfer}
For any $\lambda>0$, let $c\in C$ be an epoch initial point and set $p=P_\lambda(c)$.
Let $V_c$ denote the remainder at $x^*$ in Theorem~\ref{thm:di_bppa}, evaluated with
$w_k=c$ and $x_{k+1}=p$. Set $A_1=\lambda$, $x_1=z_1=p$, and $E_1=\lambda(F(p)-F(x^*))+D_\omega(x^*,p)$. Then
\begin{equation}\label{eq:safeguard_first_step_identity}
E_1+D_\omega(p,c)-\lambda D_f(p,c)+\lambda V_c
=D_\omega(x^*,c).
\end{equation}
Consequently,
\begin{equation}\label{eq:adaptive_epoch_initial}
E_1+(1-\lambda L)D_\omega(p,c)\le D_\omega(x^*,c).
\end{equation}
If $\lambda L\le1$, then after $\ell\ge1$ accepted calls in this epoch, its current state $(x,z,A_\ell)$ satisfies
\[
A_\ell(F(x)-F(x^*))+D_\omega(x^*,z)
\le D_\omega(x^*,c).
\]
\end{lemma}
\begin{proof}
Let $g_c=\lambda^{-1}(\nabla\omega(c)-\nabla\omega(p)) \in\nabla f(c)+\partial h(p)$. The three-point identity gives
\[
D_\omega(x^*,p)-D_\omega(x^*,c)
=\lambda\langle g_c,x^*-p\rangle-D_\omega(p,c).
\]
Combining this with $F(p)-F(x^*)+\langle g_c,x^*-p\rangle=D_f(p,c)-V_c$ proves \eqref{eq:safeguard_first_step_identity}. Relative smoothness gives $D_\omega(p,c)-\lambda D_f(p,c) \ge(1-\lambda L)D_\omega(p,c)$. Dropping $V_c\ge0$ proves \eqref{eq:adaptive_epoch_initial}. When $\lambda L\le1$, this yields $E_1\le D_\omega(x^*,c)$. Each accepted trial decreases the energy by Proposition~\ref{prop:adaptive_safeguard}, proving the endpoint estimate.
\end{proof}
\begin{remark}[Backtracking at epoch initialization]
At epoch initialization,
\eqref{eq:safeguard_first_step_identity} shows that the
computable condition
\[
\lambda D_f(p,c)\le D_\omega(p,c)
\]
is sufficient for $E_1\le D_\omega(x^*,c)$.
This condition is computable without $L$ and is guaranteed
when $\lambda L\le1$, which suggests a backtracking version of Algorithm~\ref{alg:adaptive_bppa}.
\end{remark}
Chaining Lemma~\ref{lem:adaptive_epoch_transfer} across epochs yields global oracle-complexity bounds without the trajectory-dependent conditions of Section~\ref{sec:analysis}. 
\begin{theorem}[Oracle complexity and epoch-dependent rates]
\label{thm:adaptive_bppa_complexity}
Algorithm~\ref{alg:adaptive_bppa} is well-defined for every $c_1\in C$, fixed $\lambda>0$ with $\lambda L\le1$, and integer $T\ge1$. After $T$ oracle calls, let $K_T$ be the number of nonempty epochs, including the current partial epoch when nonempty, and let $\ell_e\ge1$, $e=1,\ldots,K_T$, be their accepted-call counts. Let $A_j$ be the weights in \eqref{eq:param} with $A_1=\lambda$, and put $D_0=D_\omega(x^*,c_1)$. Then
\begin{equation}\label{eq:adaptive_posterior_rate}
0\le F(\widehat x_T)-F(x^*)
\le\frac{D_0}{\sum_{e=1}^{K_T} A_{\ell_e}},
\end{equation}
and
\begin{equation}\label{eq:adaptive_oracle_rate}
F(\widehat x_T)-F(x^*)
\le\min\left\{
\frac{2D_0}{\lambda T},\,
\frac{4D_0K_T}{\lambda T^2}
\right\}.
\end{equation}
Thus the best accepted point has a global $\cO(T^{-1})$ bound. If $\sum_e\ell_e^2=\Omega(T^2)$, in particular if the number of epochs remains bounded, the bound is $\cO(T^{-2})$.
\end{theorem}
\begin{proof}
By Lemma~\ref{lem:prox_interior} and $T_\tau(C\times C)\subset C$, the algorithm is well-defined with finite trial margins, and its first oracle output is accepted by construction.

For epoch $e$, let $c_e$ be its initial point and $(x^{(e)},z^{(e)})$ its last accepted state. Lemma~\ref{lem:adaptive_epoch_transfer} gives
\[
A_{\ell_e}(F(x^{(e)})-F(x^*))+D_\omega(x^*,z^{(e)})
\le D_\omega(x^*,c_e).
\]
Since $c_{e+1}=z^{(e)}$, summation over the nonempty epochs yields
\begin{equation}\label{eq:adaptive_epoch_telescope}
\sum_{e=1}^{K_T} A_{\ell_e}(F(x^{(e)})-F(x^*))
+D_\omega(x^*,z^{(K_T)})\le D_0.
\end{equation}
Since $F(\widehat x_T)\le F(x^{(e)})$ for every $e$, we have \eqref{eq:adaptive_posterior_rate}.

Each nonempty epoch contains at least one accepted call and at most one rejected call, including a rejection on the final call. Thus
\[
T\le\sum_{e=1}^{K_T}\ell_e+K_T
\le2\sum_{e=1}^{K_T}\ell_e.
\]
By \eqref{eq:weight_lower_bound}, $A_j\ge\lambda(j+1)^2/4\ge\lambda j$. 
Hence by Cauchy--Schwarz, we have
\[
\begin{aligned}
\sum_{e=1}^{K_T}A_{\ell_e}
&\ge\lambda\sum_{e=1}^{K_T}\ell_e\ge\frac{\lambda T}{2},\\
\sum_{e=1}^{K_T}A_{\ell_e}
&\ge\frac{\lambda}{4}\sum_{e=1}^{K_T}(\ell_e+1)^2
\ge\frac{\lambda}{4K_T}
\left(\sum_{e=1}^{K_T}\ell_e+K_T\right)^2
\ge\frac{\lambda T^2}{4K_T}.
\end{aligned}
\]
Substituting these bounds into \eqref{eq:adaptive_posterior_rate} proves \eqref{eq:adaptive_oracle_rate}. The quadratic weight bound also gives $\cO(T^{-2})$ whenever $\sum_e\ell_e^2=\Omega(T^2)$.
\end{proof}

\subsection{Finite Restart and Eventual Linear Convergence}
For the finite-restart and local convergence results, we assume quadratic growth at an interior minimizer, which is weaker than local strong convexity and has been used to establish linear convergence~\citep{drusvyatskiy2018error,fercoq2019adaptive}.
\begin{assumption}[Local quadratic growth]
\label{ass:local_qg}
There exist $m_F>0$ and $r>0$ with $B_r(x^*)\subset C$ such that
\begin{equation}\label{eq:local_qg}
F(x)-F(x^*)\ge\frac{m_F}{2}\|x-x^*\|^2,
\qquad x\in B_r(x^*)\cap\dom F.
\end{equation}
\end{assumption}

The local margin bound in Appendix~\ref{app:local_analysis} ensures acceptance from small initial energy. First-step dissipation then bounds the total number of restarts from any initial point in $C$.

\begin{corollary}[Local acceptance and linear convergence]\label{lem:di_local_acceptance}
Under Assumption~\ref{ass:local_qg}, suppose that $\lambda L<1$ and that $\omega$ is $C^2$ near $x^*$ with $\nabla^2\omega(x^*)\succ0$.
There exist $\varepsilon_*>0$ and $\rho\in(0,1)$ such that any epoch of Algorithm~\ref{alg:adaptive_bppa}, initialized with $x_1=z_1=p$ and $A_1=\lambda$, satisfying
\[
E_1=\lambda(F(p)-F(x^*))+D_\omega(x^*,p)\le\varepsilon_*,
\]
accepts every subsequent trial and satisfies, for every epoch index $j\ge1$,
\[
E_j\le E_1\rho^{j-1},\qquad
0\le F(x_j)-F(x^*)\le\frac{E_1}{A_j}\rho^{j-1}.
\]
In particular, $D_\omega(x^*,c)\le\varepsilon_*$ at the epoch's initial point $c$ is sufficient.
\end{corollary}
\begin{proof}
By Theorem~\ref{thm:local_bppa_stability} and
Corollary~\ref{cor:local_bppa_linear}, the unrestarted DI sequence
from the same initial state satisfies
$E_j\le E_1\rho^{j-1}\le\varepsilon_*$.
At each of these states, \eqref{eq:local_margin_control} gives
$\mathcal M_j=\widehat Q_j+U_j-J_j\ge0$.
Thus S-DI accepts each trial by induction and inherits the stated
linear bounds. The sufficient condition on $c$ follows from
$E_1\le D_\omega(x^*,c)$ in \eqref{eq:adaptive_epoch_initial}.
\end{proof}
\begin{theorem}[Finite restart under local quadratic growth]
\label{thm:safeguard_finite_restart}
Under standing assumptions and Assumption~\ref{ass:local_qg}, suppose that $\lambda L<1$ and that $\omega$ is $C^2$ near $x^*$ with $\nabla^2\omega(x^*)\succ0$. Let $c_1\in C$, $D_0=D_\omega(x^*,c_1)$, and take $\varepsilon_*$ from Corollary~\ref{lem:di_local_acceptance}. 
When $D_0>\varepsilon_*$, define $\mathcal K_0:=\{c\in C:\varepsilon_*\le D_\omega(x^*,c)\le D_0\}$, then
\[
d_{\min}:=\min_{c\in\mathcal K_0}D_\omega(P_\lambda(c),c)>0.
\]
And Algorithm~\ref{alg:adaptive_bppa} has a total number $R_*$ of rejected trials satisfying $R_*\le\overline R$, where
\begin{equation}\label{eq:restart_count_bound}
\overline R:=
\begin{cases}
0,&D_0\le\varepsilon_*,\\[1mm]
\displaystyle\left\lfloor\frac{D_0}{(1-\lambda L)d_{\min}}\right\rfloor,
&D_0>\varepsilon_*.
\end{cases}
\end{equation}
Consequently, for every $T\ge1$, its best accepted point satisfies
\begin{equation}\label{eq:safeguard_finite_restart_rate}
0\le F(\widehat x_T)-F(x^*)
\le\frac{4D_0(\overline R+1)}{\lambda T^2}.
\end{equation}
\end{theorem}
\begin{proof}
Set $\sigma=1-\lambda L>0$. For epoch $e$, let $c_e$ be its initial point, $p_e=P_\lambda(c_e)$, and $D_e=D_\omega(x^*,c_e)$.
By \eqref{eq:safeguard_first_step_identity}, relative smoothness, and the decrease of accepted energies, every completed epoch satisfies
\begin{equation}\label{eq:restart_anchor_decrease}
D_{e+1}\le E_{1,e}
\le D_e-\sigma D_\omega(p_e,c_e)-\lambda V_{c_e},
\end{equation}
where $E_{1,e}$ is the initial energy of epoch $e$ and $V_{c_e}\ge0$. Thus $D_e\le D_0$. If $D_e\le\varepsilon_*$, Corollary~\ref{lem:di_local_acceptance} rules out any rejection in that epoch. This proves $R_*=0$ when $D_0\le\varepsilon_*$.

Suppose $D_0>\varepsilon_*$. Lemma~\ref{lem:bregman_sublevel_compact} makes $\mathcal K_0$ compact, and $c_1\in\mathcal K_0$. 
By Lemma~\ref{lem:prox_interior} applied with $\varphi=h$ and the parameter $t(c)=\nabla\omega(c)-\lambda\nabla f(c)$, which is continuous on $C$, the map $P_\lambda$ is continuous on $C$.
Hence the minimum defining $d_{\min}$ is attained. If it were zero, some $c\in\mathcal K_0$ would satisfy $P_\lambda(c)=c$, so optimality gives $-\nabla f(c)\in\partial h(c)$. Thus $c$ is optimal. Convexity and Assumption~\ref{ass:local_qg} make $x^*$ the unique minimizer, contradicting $D_\omega(x^*,c)\ge\varepsilon_*>0$. Therefore $d_{\min}>0$.

Every epoch $e$ that ends in a rejection has $D_e>\varepsilon_*$, so $c_e\in\mathcal K_0$ and \eqref{eq:restart_anchor_decrease} yields $D_e-D_{e+1}\ge\sigma d_{\min}$. Summing over any finite number $r$ of rejected epochs gives
\[
r\sigma d_{\min}\le D_0-D_{r+1}\le D_0.
\]
This proves $R_*\le\overline R<\infty$. Finally, $K_T\le R_*+1\le\overline R+1$ in \eqref{eq:adaptive_oracle_rate} proves \eqref{eq:safeguard_finite_restart_rate}.
\end{proof}

It can be shown that the  quadratically regularized Burg kernel satisfies the local kernel regularity requirement even though it does not satisfy
global quadratic triangle scaling; see
Remark~\ref{app:burg_scaling_separation}.

The final epoch of S-DI need not start with small energy.
We next show that its energy eventually falls below the threshold required for local linear convergence.
\begin{corollary}[Eventual linear convergence of the safeguarded DI methods]
\label{cor:bppa_eventual_linear}
Under the assumptions of Theorem~\ref{thm:safeguard_finite_restart},
let $j\ge1$ index the states in the final epoch of
Algorithm~\ref{alg:adaptive_bppa}.
There exist $j_0\ge1$ and $\rho\in(0,1)$ such that, for every $j\ge j_0$,
\begin{equation}\label{eq:bppa_eventual_linear}
\begin{aligned}
E_j&\le E_{j_0}\rho^{j-j_0},\\
0\le F(x_j)-F(x^*)
&\le\frac{E_{j_0}}{A_j}\rho^{j-j_0}
\le\frac{4E_{j_0}\rho^{j-j_0}}{\lambda(j+1)^2}.
\end{aligned}
\end{equation}
\end{corollary}
\begin{proof}
    The proof is given in Appendix~\ref{app:eventual_linear_proof}.
\end{proof}

\subsection{Finite Restart under Weak Sharpness}
We next use local weak sharpness, a classical condition for
finite proximal termination~\citep{burke1993weak,ferris1991finite},
to establish finite restart.

\begin{proposition}[Finite restart under local weak sharpness]
\label{prop:safeguard_weak_sharpness}
Let $\lambda>0$ satisfy $\lambda L<1$. Suppose
\[
S:=\{x\in\dom F\cap C:F(x)=F(x^*)\}\ne\varnothing.
\]
For every $q\in S$, assume there are an open neighborhood $\mathcal O_q\subset C$ of $q$ and $\alpha_q>0$ such that
\[
F(x)-F(x^*)\ge\alpha_q\operatorname{dist}(x,S),
\qquad x\in\mathcal O_q\cap\dom F.
\]
Then any exact DI method implemented by Algorithm~\ref{alg:adaptive_bppa}, run indefinitely from any $c_1\in C$, rejects only finitely many trials. All sufficiently late accepted main iterates belong to $S$.
\end{proposition}

\begin{proof}

    For $x\in(\mathcal O_q\cap\dom F)\setminus S$ and $v\in\nabla f(x)+\partial h(x)$, taking the infimum over $y\in S$ in the subgradient inequality, we have
\[
0<F(x)-F(x^*)
\le\|v\|\operatorname{dist}(x,S).
\]
Thus $\operatorname{dist}(x,S)>0$, and local weak sharpness implies $\|v\|\ge\alpha_q$.

Now suppose that $y_j$ lies in a compact subset of $C$, with $v_j\in\nabla f(y_j)+\partial h(y_j)$ and $v_j\to0$. If infinitely many $y_j$ were nonoptimal, a subsequence of these points would converge to some $q\in C$. Continuity of $\nabla f$ and the closed graph of $\partial h$ give $-\nabla f(q)\in\partial h(q)$, hence $q\in S$. For sufficiently large indices of this subsequence, the preceding bound requires $\|v_j\|\ge\alpha_q$, contradicting $v_j\to0$. 
So $y_j\in S$ for all sufficiently large $j$.

Fix $q_0\in S$ and suppose infinitely many trials are rejected. By \eqref{eq:restart_anchor_decrease}, applied with comparator $q_0$, both $D_\omega(q_0,c_e)$ and $D_\omega(q_0,p_e)$ remain bounded, while $D_\omega(p_e,c_e)\to0$. Lemma~\ref{lem:bregman_sublevel_compact} places $c_e,p_e$ in a common compact subset of $C$. Continuity and strict positivity of the divergence on this set imply $p_e-c_e\to0$.

Define
\[
v_e:=\lambda^{-1}
\bigl(\nabla\omega(c_e)-\nabla\omega(p_e)\bigr)
+\nabla f(p_e)-\nabla f(c_e).
\]
By optimality, $v_e\in\nabla f(p_e)+\partial h(p_e)$ and gradient continuity on the compact set ensures $v_e\to0$. 
The preceding identification argument gives $p_e\in S$ for all sufficiently large $e$. Starting an epoch from $x=z=p_e\in S$, uniqueness of the exact oracle keeps every subsequent update at $p_e$, with $g=0$ and $\mathcal M=0$. This epoch cannot reject a trial, contradicting infinitely many rejections.

In the final epoch, the energy bound with $q_0$ and Lemma~\ref{lem:bregman_sublevel_compact} bound $s_k=\nabla\omega(z_k)$. The reconstruction recursion $u_{k+1}=(1-\theta_k)u_k+\theta_k s_{k+1}$ also bounds $u_k$ and $\nabla\omega(w_k)$. By continuity of the inverse mirror map, $x_k,z_k,w_k$ lie in a common compact subset of~$C$. By \eqref{eq:auxiliary_dual_update} and \eqref{eq:weight_increment}, we obtain
\[
g_{k+1}=\frac{s_k-s_{k+1}}{a_{k+1}}\to0,
\qquad a_{k+1}\to\infty.
\]
From \eqref{eq:oracle_residual}, $\nabla\omega(w_k)-u_{k+1}=\lambda g_{k+1}\to0$. Uniform continuity of $\nabla H$ on the corresponding dual compact set implies $w_k-x_{k+1}\to0$. Therefore
\[
v_{k+1}:=g_{k+1}
+\nabla f(x_{k+1})-\nabla f(w_k)
\in\nabla f(x_{k+1})+\partial h(x_{k+1}),
\qquad v_{k+1}\to0.
\]
Applying the same identification argument to $(x_{k+1},v_{k+1})$ shows that all sufficiently late main iterates belong to $S$. The incumbent is therefore optimal after finitely many calls.
\end{proof}

Note that this proposition requires neither local quadratic growth nor $C^2$ regularity of the kernel. 

\section{Non-Restarting Safeguard with Dual Substeps}
\label{ssec:fractional_fallback}

S-DI resets the weights after a rejected trial. DI-FS as proposed in Algorithm~\ref{alg:di_fs}  instead keeps~\eqref{eq:param} and uses the same proximal oracle and acceptance margin on a sequence of dual substeps. Throughout this section, we assume $\lambda L\le1$.

\subsection{Dual Substeps}
Fix an outer step with $A=A_k$, $A^+=A_{k+1}$, $\theta=\theta_k$, and $a=A^+-A$. At a fallback state $(q,\zeta)$, a fraction $\alpha\in[\theta,1]$ defines the interpolation parameter and weights
\begin{equation}\label{eq:dual_substep_weights}
\vartheta=\frac{\theta}{\alpha},\qquad
B^+=\frac{\lambda}{\vartheta^2},\qquad
B=(1-\vartheta)B^+,\qquad
b=B^+-B=\alpha a.
\end{equation}
Use the DI trial $w=T_\vartheta(q,\zeta)$, $p=P_\lambda(w)$, and $\zeta^+=T_{1/\vartheta}(q,p)$. The substep margin $\mathcal M_B$ is~\eqref{eq:adaptive_margin} evaluated at the current state $(B,q,\zeta)$ and the trial state $(B^+,p,\zeta^+)$, with $b$ in place of $a$ and the residual from~\eqref{eq:oracle_residual}. Thus the proximal parameter remains $\lambda$; only the interpolation parameter and the weights in the test change. 
For $B>0$, the objective comparison after acceptance uses the values already needed for the margin.

\begin{algorithm}[htbp]
\caption{Dual-inertial Bregman algorithm with fractional substeps (DI-FS)}
\label{alg:di_fs}
\begin{algorithmic}[1]
\REQUIRE $c_1\in C$, $\lambda>0$ with $\lambda L\le1$, and an outer horizon $N\ge1$.
\STATE Compute $p=P_\lambda(c_1)$ and set $(x_1,z_1,A_1)=(p,p,\lambda)$.
\FOR{$k=1,\ldots,N-1$}
    \STATE Compute $A_{k+1}$ and $\theta_k$ from~\eqref{eq:param}.
    \STATE Generate the DI trial $(w_k,x^+,z^+)$ and compute $\mathcal M$ by~\eqref{eq:adaptive_margin}.
    \IF{$\mathcal M\ge0$}
        \STATE Set $(x_{k+1},z_{k+1})=(x^+,z^+)$.
    \ELSE
        \STATE Set $(q,\zeta)=(x_k,z_k)$, $r=1$, and $\alpha=\max\{1/2,\theta_k\}$.
        \WHILE{$r>0$}
            \STATE Compute~\eqref{eq:dual_substep_weights}, the trial $(w,p,\zeta^+)$, and $\mathcal M_B$.
            \IF{$\mathcal M_B\ge0$}
                \STATE Replace $q$ by $p$ if $F(p)<F(q)$.
                \STATE Set $\zeta=\zeta^+$ and $r=r-\alpha$; keep $\alpha$ unchanged.
            \ELSE
                \STATE Set $\alpha=\max\{\theta_k,\alpha/2\}$.
            \ENDIF
        \ENDWHILE
        \STATE Set $(x_{k+1},z_{k+1})=(q,\zeta)$.
    \ENDIF
\ENDFOR
\RETURN $x_N$.
\end{algorithmic}
\end{algorithm}

The rejected outer trial corresponds to $\alpha=1$ and is not repeated. The last fraction is not truncated to $r$: the loop stops once the accepted fractions sum to at least one, while the outer weight remains $A^+$.

\begin{lemma}[Energy decrease of the fallback]\label{lem:frac_identity}
Every completed fallback in Algorithm~\ref{alg:di_fs} satisfies $E_{k+1}\le E_k$ for the energy~\eqref{eq:lyap_bppa}.
\end{lemma}
\begin{proof}
Write $\phi=F-F(x^*)$ and index the accepted substeps and their fractions $\alpha_j$ by $j=0,\ldots,S-1$, with $(q_0,\zeta_0)=(x_k,z_k)$. Denote the trial main point by $p_j$, the margin by $\mathcal M_j$, and set $D_j=D_\omega(x^*,\zeta_j)$. If $B_j>0$, apply~\eqref{eq:adaptive_edge_energy} with weights $B_j,B_j^+$. If $B_j=0$, then $\vartheta=1$, $w=\zeta_j$, and $\zeta_{j+1}=p_j$; apply the first-step identity~\eqref{eq:safeguard_first_step_identity}. Dropping the nonnegative remainder gives in both cases
\[
B_j^+\phi(p_j)+D_{j+1}+\mathcal M_j
\le B_j\phi(q_j)+D_j.
\]
Since $\phi(q_{j+1})=\min\{\phi(q_j),\phi(p_j)\}$ and
$B_j=\alpha_j(\alpha_j-\theta)A^+\le(1-\theta)A^+=A$, it follows that
\[
\begin{aligned}
b_j\phi(q_{j+1})+D_{j+1}-D_j+\mathcal M_j
\le B_j\bigl[\phi(q_j)-\phi(q_{j+1})\bigr]
\le A\bigl[\phi(q_j)-\phi(q_{j+1})\bigr].
\end{aligned}
\]
Summing over $j=0,\ldots,S-1$ telescopes $D_{j+1}-D_j$ and $\phi(q_j)-\phi(q_{j+1})$. Also, $\phi(q_{j+1})\ge\phi(q_S)\ge0$ and $\sum_j b_j=a\sum_j\alpha_j\ge a$. Thus, dropping $\sum_j\mathcal M_j\ge0$ yields $(A+a)\phi(q_S)+D_S\le A\phi(q_0)+D_0$, which is $E_{k+1}\le E_k$.
\end{proof}

\subsection{Termination and Oracle Complexity}

Let $m_k$ count all proximal calls at outer update $k$, including the initial DI trial and any rejected substep trials. Through outer state $N$, including initialization, the total and additional fallback counts are
\[
T_N=1+\sum_{k=1}^{N-1}m_k=N+C_N,
\qquad C_N=\sum_{k=1}^{N-1}(m_k-1).
\]

\begin{theorem}[Global rates of DI-FS]\label{thm:frac_hybrid}
Consider Algorithm~\ref{alg:di_fs} with $c_1\in C$ and $\lambda>0$ satisfying $\lambda L\le1$.
\begin{enumerate}
\item Every fallback terminates. For every $N\ge1$, $T_N\le N^2$ and
\begin{equation}\label{eq:frac_hybrid_rate}
F(x_N)-F(x^*)\le\frac{E_1}{A_N}
\le\frac{4E_1}{\lambda(N+1)^2}\le\frac{4E_1}{\lambda T_N},
\qquad E_1\le D_\omega(x^*,c_1).
\end{equation}
\item If $C_k\le\bar c k$ for all $k\ge1$ and some $\bar c\ge0$, then for every $N\ge1$,
\[
F(x_N)-F(x^*)\le\frac{4(1+\bar c)^2E_1}{\lambda T_N^2}.
\]
\end{enumerate}
Consequently, given any integer budget of $T\ge1$ proximal calls, stopping before making a call that would exceed the budget returns the last completed outer iterate $x_N$ ($N\ge1$) satisfying
\[
F(x_N)-F(x^*)\le\frac{4E_1}{\lambda T} \quad \text{or} \quad
F(x_N)-F(x^*)\le\frac{4(1+\bar c)^2E_1}{\lambda T^2}.
\]
\end{theorem}

\begin{proof}
At $\alpha=\theta_k$, \eqref{eq:dual_substep_weights} gives $\vartheta=1$, $B=0$, and $b=\lambda$. Substituting $w=\zeta$ and $\zeta^+=p$ into the margin gives
\[
\mathcal M_0=D_\omega(p,\zeta)-\lambda D_f(p,\zeta)
\ge(1-\lambda L)D_\omega(p,\zeta)\ge0.
\]
Counting the initial DI trial as $\alpha=1$, the fraction before the $r$-th rejection is $\max\{\theta_k,2^{-(r-1)}\}$, since only rejections reduce $\alpha$. Trials at $\alpha=\theta_k$ are always accepted, so a rejection requires $r-1<\log_2(1/\theta_k)$. Hence outer step $k$ has at most $\lceil\log_2(1/\theta_k)\rceil$ rejected trials, including the initial DI trial.
Each accepted substep uses a fraction of at least $\theta_k$, so the loop ends after at most $\lceil1/\theta_k\rceil$ acceptances. Since $\theta_k^{-1}=t_{k+1}\le k+1$ by~\eqref{eq:param} with $t_1=1$, we have $m_k\le2k+1$ whether or not fallback occurs. Hence $T_N\le1+\sum_{k=1}^{N-1}(2k+1)=N^2$.

The energy does not increase in an accepted outer DI trial by Proposition~\ref{prop:adaptive_safeguard}, or in a fallback by Lemma~\ref{lem:frac_identity}. Thus $E_N\le E_1$, and \eqref{eq:weight_lower_bound} gives \eqref{eq:frac_hybrid_rate}, while the initial energy bound follows from~\eqref{eq:adaptive_epoch_initial}. Given a budget of $T\ge1$ proximal calls, initialization uses $T_1=1\le T$ calls, so at least one outer iterate is completed. The call counts of the uncapped run satisfy $T_N\le T<T_{N+1}\le(N+1)^2$, and~\eqref{eq:frac_hybrid_rate} yields $F(x_N)-F(x^*)\le\frac{4E_1}{\lambda T}$.

For part~(2), substituting $T_N=N+C_N\le(1+\bar c)N<(1+\bar c)(N+1)$ into the outer-iteration bound in~\eqref{eq:frac_hybrid_rate} gives the $T_N^{-2}$ bound. Under the budget $T\ge1$, the same condition at $k=N+1$ gives $T<T_{N+1}\le(1+\bar c)(N+1)$, and~\eqref{eq:frac_hybrid_rate} yields the $T^{-2}$ bound.
\end{proof}

Note that Theorem~\ref{thm:frac_hybrid} requires neither triangle scaling nor convexity of $D_\omega$ in its second argument or joint convexity, and it permits $x^*\in\partial\dom\omega$. 
\section{Numerical Experiments}\label{sec:experiments_numeric}
\label{sec:numerical_experiments}
We compare DI-BPG and its S-DI and DI-FS variants on relative-entropy regression, and DI-BPPA and its S-DI and DI-FS variants on nonnegative least squares. Each comparison uses the same instances and initial points. Experiments were run on a Mac mini (Apple M4) with one numerical-library thread, using Python 3.12.14. Timings are pointwise medians of three serial runs and include function evaluations, gain searches, Newton solves, and safeguard tests. Throughout the experiments, performance is measured by the relative objective gap where objective gaps are normalized by the gap at the common initial point.

Both problems use the orthant entropy $\omega(x)=\sum_i(x_i\log x_i-x_i)$ on $\R_+^n$, with $0\log0=0$, and its Kullback--Leibler divergence
\[
D_{\rm KL}(v,w)
=\sum_i\left(v_i\log\frac{v_i}{w_i}-v_i+w_i\right),
\qquad v\ge0,\quad w>0.
\]
For $H(u)=\sum_i e^{u_i}$, the directional third derivative along a dual segment has nonnegative derivative $\sum_i e^{u_{k-1,i}+td_{k,i}}d_{k,i}^4$. Thus the supremum in Theorem~\ref{prop:directional_budget} is attained at $t=1$:
\[
\eta_k=\left[\sum_i z_{k,i}d_{k,i}^3\right]_+,
\qquad d_k=s_k-u_{k-1}.
\]
For seed 1, we plot $\eta_k$, $\sum_{k=2}^N\theta_{k-1}\eta_k$, and $\mathcal B_N$ for DI methods, together with cumulative S-DI rejections, for each experiment. 

\subsection{Relative-Entropy Regression}
\begin{figure}[ht]
\centering
\includegraphics[width=0.9\textwidth]{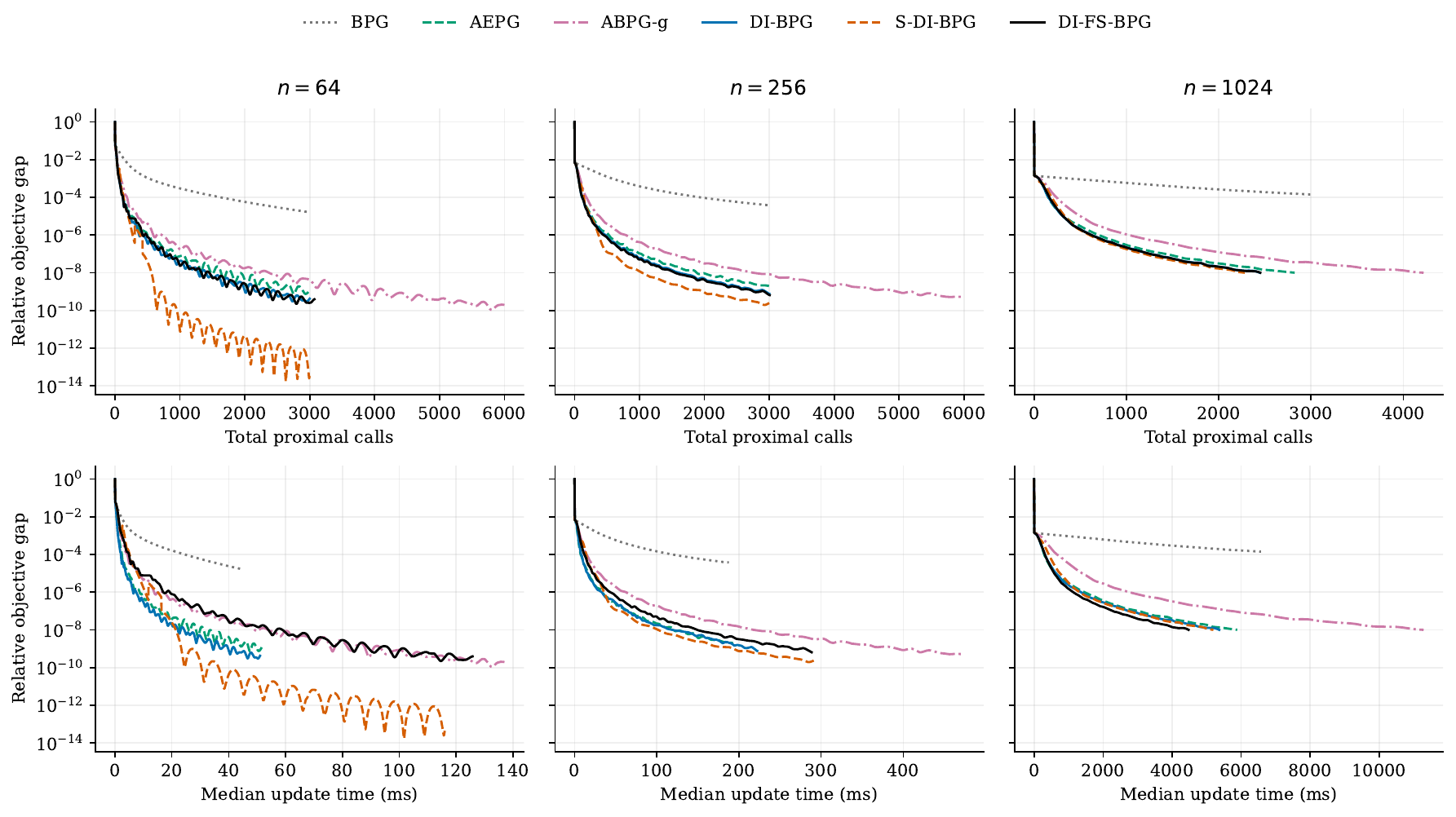}
\caption{Relative-entropy regression, seed 1 with $n=64,256, 1024$.}
\label{fig:orthant_comparison}
\end{figure}

Following the problem class in Hanzely et al.~\citep[Section~6.3]{hanzely2021accelerated}, we minimize
\[
F(x)=D_{\rm KL}(Ax,b)+\tau\mathbf1^\top x,\qquad x\ge0,
\quad \tau=0.1.
\]
Both the seed-1 plots and the multi-instance comparison use $n=64,256,1024$, with $m=3n$ in all cases. The entries of $A$ are drawn uniformly on $(0,1)$ and its columns are normalized to sum to one, giving relative smoothness with $L=1$. Independently, $\log x_i^*$ is uniform on $[\log0.1,\log3]$, and $b=e^\tau Ax^*$, so $\nabla F(x^*)=0$ and
\[
F(x)-F(x^*)=D_{\rm KL}(Ax,Ax^*).
\]
The seed-1 comparison includes BPG, AEPG~\citep{merchav2026adaptive}, ABPG-g~\citep{hanzely2021accelerated}, DI-BPG, S-DI-BPG, and DI-FS-BPG, starting at $0.5\mathbf1$ for 3000 outer iterations  (stopping early at the $10^{-8}$ target for $n=1024$ in Figure~\ref{fig:orthant_comparison}). BPG, AEPG, DI-BPG, S-DI-BPG, and DI-FS-BPG use $\lambda=0.99/L$. AEPG retains the moving constraint $x^+\ge(1-\theta)x$. ABPG-g uses Algorithm~3 of the paper, with $1/L$, exponent two, initial gain one, gain multiplier 1.5, and minimum gain $10^{-3}$. DI-FS-BPG follows Algorithm~\ref{alg:di_fs}; its call count includes rejected DI trials and all fractional-substep trials.

\begin{figure}[ht]
\centering
\includegraphics[width=0.9\textwidth]{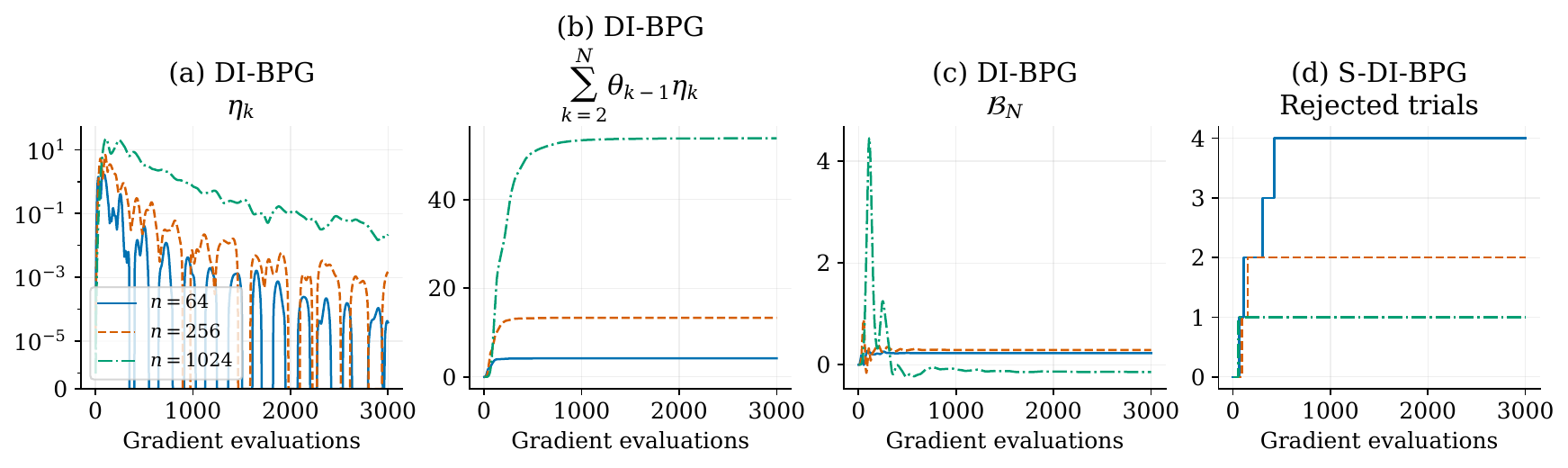}
\caption{Relative-entropy regression, seed 1: geometry of unrestarted DI-BPG in (a)--(c) and S-DI-BPG rejections in (d), for $n=64,256,1024$}
\label{fig:regression_geometry}
\end{figure}

In Figure~\ref{fig:orthant_comparison}, DI-BPG reached smaller final gaps than AEPG at comparable update time. ABPG-g attained a slightly smaller final objective gap than unrestarted DI-BPG, but used almost twice as many proximal calls and more time. S-DI reduced the final objective gaps, at a higher computational cost than DI-BPG. DI-FS-BPG required only a small number of fallback substeps and attained final objective gaps comparable to DI-BPG.
Over the plotted iteration ranges in
Figure~\ref{fig:regression_geometry}, the weighted sums
$\sum_{k=2}^N\theta_{k-1}\eta_k$ and connection budgets
$\mathcal B_N$ showed little late-run variation in all three sizes.

Table~\ref{tab:regression_target_times} supplements the seed-1 comparison with 20 random instances at each size. For each instance, we use the median of three serial timings; the table reports medians across instances. Within the 3000 outer-iteration limit, standard BPG failed to reach either target ($10^{-6}$ or $10^{-8}$) on any of the 20 instances at any problem size. Although ABPG-g used the fewest outer iterations, DI-BPG has shorter median target times than AEPG and ABPG-g at every size. Because each gain search may evaluate multiple candidate steps, ABPG-g required roughly twice as many total proximal calls as outer iterations, whereas S-DI-BPG achieved the fewest median total proximal calls at every problem size and target tolerance. Moreover, DI-BPG attained the shortest median target times for $n=64,256$, while DI-FS-BPG attained the shortest median target times at both targets for $n=1024$.  

\begin{table}[htbp]
\centering
\caption{Median target times, outer iterations, and total proximal calls over 20 random instances}
\label{tab:regression_target_times}
\scriptsize
\renewcommand{\arraystretch}{0.95}
\setlength{\tabcolsep}{1.2pt}
\begin{tabular*}{\textwidth}{@{\extracolsep{\fill}}llrrrrrrrrr@{}}
\toprule
 & & \multicolumn{3}{c}{$n=64$} & \multicolumn{3}{c}{$n=256$} & \multicolumn{3}{c}{$n=1024$} \\
\cmidrule(lr){3-5}\cmidrule(lr){6-8}\cmidrule(lr){9-11}
Method & Target & \shortstack{Time\\(ms)} & \shortstack{Outer\\iters.} & \shortstack{Total\\calls} & \shortstack{Time\\(ms)} & \shortstack{Outer\\iters.} & \shortstack{Total\\calls} & \shortstack{Time\\(ms)} & \shortstack{Outer\\iters.} & \shortstack{Total\\calls} \\
\midrule
BPG & $10^{-6}$ & --- & --- & --- & --- & --- & --- & --- & --- & --- \\
BPG & $10^{-8}$ & --- & --- & --- & --- & --- & --- & --- & --- & --- \\
\addlinespace[1pt]
AEPG & $10^{-6}$ & 5.5 & 346.5 & 346.5 & 31.6 & 491.5 & 491.5 & 1404.1 & 687.5 & 687.5 \\
AEPG & $10^{-8}$ & 18.2 & 1111 & 1111 & 117.8 & 1817.5 & 1817.5 & 5602.6 & 2783 & 2783 \\
\addlinespace[1pt]
ABPG-g & $10^{-6}$ & 11.1 & \textbf{249} & 497 & 58.8 & \textbf{378} & 756 & 2481.0 & \textbf{514} & 1027.5 \\
ABPG-g & $10^{-8}$ & 37.0 & \textbf{825.5} & 1650 & 218.9 & \textbf{1380} & 2759.5 & 10136.7 & \textbf{2116} & 4231.5 \\
\addlinespace[1pt]
DI-BPG & $10^{-6}$ & \textbf{4.8} & 294 & 294 & \textbf{29.3} & 448 & 448 & 1251.8 & 627.5 & 627.5 \\
DI-BPG & $10^{-8}$ & \textbf{15.4} & 949.5 & 949.5 & \textbf{103.1} & 1603.5 & 1603.5 & 5070.7 & 2495 & 2495 \\
\addlinespace[1pt]
S-DI-BPG & $10^{-6}$ & 10.8 & 274.5 & \textbf{274.5} & 39.8 & 410.5 & \textbf{410.5} & 1309.0 & 598.5 & \textbf{598.5} \\
S-DI-BPG & $10^{-8}$ & 32.8 & 853 & \textbf{853} & 134.3 & 1380.5 & \textbf{1380.5} & 4223.0 & 2282.5 & \textbf{2282.5} \\
\addlinespace[1pt]
DI-FS-BPG & $10^{-6}$ & 10.4 & 290.5 & 293.5 & 42.7 & 441.5 & 445.5 & \textbf{1042.0} & 611 & 624 \\
DI-FS-BPG & $10^{-8}$ & 33.1 & 943 & 946 & 146.2 & 1542.5 & 1545.5 & \textbf{4137.7} & 2443 & 2456 \\
\bottomrule
\end{tabular*}
\par\smallskip
\begin{minipage}{\linewidth}\footnotesize
Proximal calls (``Total calls'') include rejected outer and substep trials, and gain-search trials; outer iterations (``Outer iters.'') exclude inner gain-search trials.
\end{minipage}
\end{table}

\subsection{Nonnegative Least Squares}
We consider nonnegative least squares as in~\citep{lawson1995least,slawski2013nonnegative},
\[
\min_{x\ge0} F(x)=\tfrac12\|Ax-b\|^2,
\qquad A\in\R^{480\times320}.
\]
Before column normalization, the design satisfies $A_{:,1}=q_1$ and $A_{:,j}=\rho A_{:,j-1}+\sqrt{1-\rho^2}\,q_j$, where $q_j$ are independent standard Gaussian vectors. We use $\rho=0.3,0.95$ and $b=Ax^{\rm true}+0.01\varepsilon$, with independent standard Gaussian $\varepsilon$. The support of $x^{\rm true}$ is a uniformly chosen set of 80 coordinates, whose values are uniform on $[0.2,1.2]$. All methods start at $0.2\mathbf1$, use the fixed parameter $\lambda=1$, and count 500 proximal calls, including rejections. We compare BPPA, V-iBPPA, DI-BPPA, S-DI-BPPA, and DI-FS-BPPA. V-iBPPA follows Algorithm~2 of Yang and Toh~\citep{yang2022bregman}, using exponent two and $c_0=\pi=\gamma=\tau_1=1$ initially. Following their Remark~5.2, $\tau_1$ is doubled whenever $\tau_1\theta_k<0.1$, and $\theta_k$ is recomputed. All methods approximate the proximal oracle using the same damped Newton solver, with infinity-norm residual tolerance $10^{-12}$ for the logarithmic optimality equation
$
u-\log w+\lambda A^\top(Ae^u-b)=0$ and $x^+=e^u.
$
An independent SciPy solve provides $x^*$, with projected-gradient KKT residual below $2\times10^{-14}$. The reference minimizers have zero coordinates and therefore lie outside $C=\mathbb R_{++}^n$. The global bounds in Theorems~\ref{thm:adaptive_bppa_complexity} and~\ref{thm:frac_hybrid} allow such boundary minimizers, whereas the local convergence results require $x^*\in C$.
Objective gaps are evaluated as
\[
F(x)-F(x^*)=\tfrac12\|A(x-x^*)\|^2
+\langle A^\top(Ax^*-b),x-x^*\rangle.
\]
The single-step balance \eqref{eq:adaptive_edge_energy} and its cumulative epoch balance had normalized residuals below $4\times10^{-11}$ in all DI runs, using the maximum of one and the sum of constituent magnitudes as the normalization.

\begin{figure}[htbp]
\centering
\includegraphics[width=0.9\textwidth]{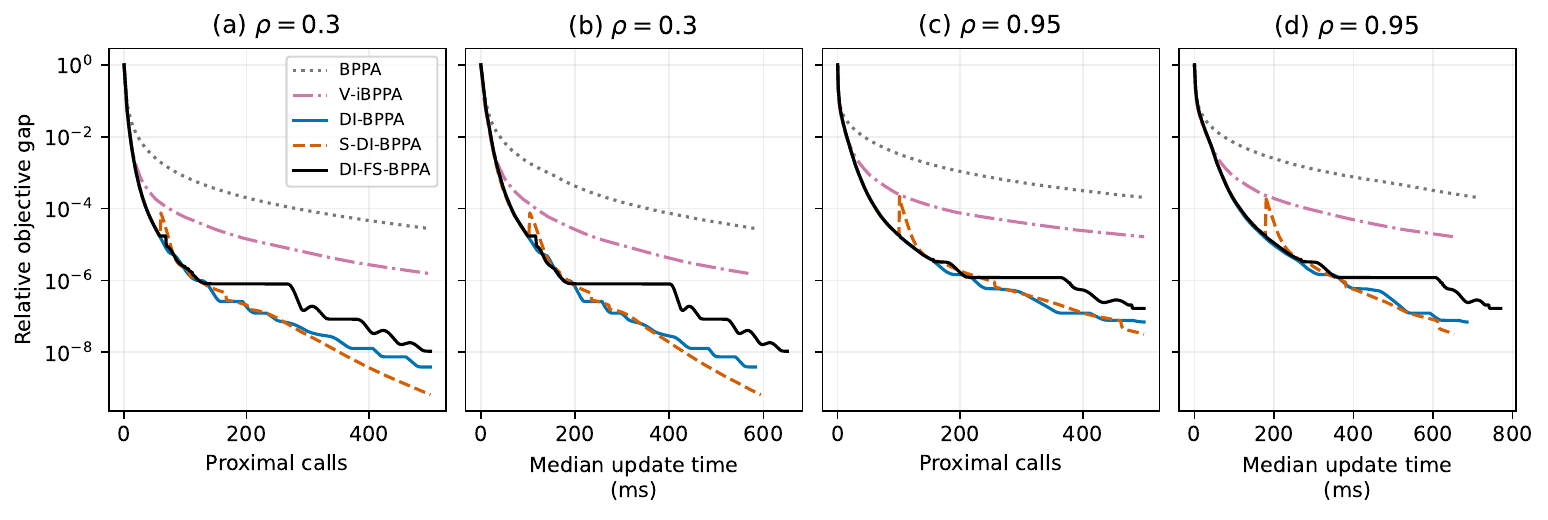}
\caption{Nonnegative least squares, $480\times320$, seed 1: relative objective gaps against total proximal calls and median cumulative update time. 
}
\label{fig:nnls_comparison}
\end{figure}

\begin{figure}[htbp]
\centering
\includegraphics[width=0.9\textwidth]{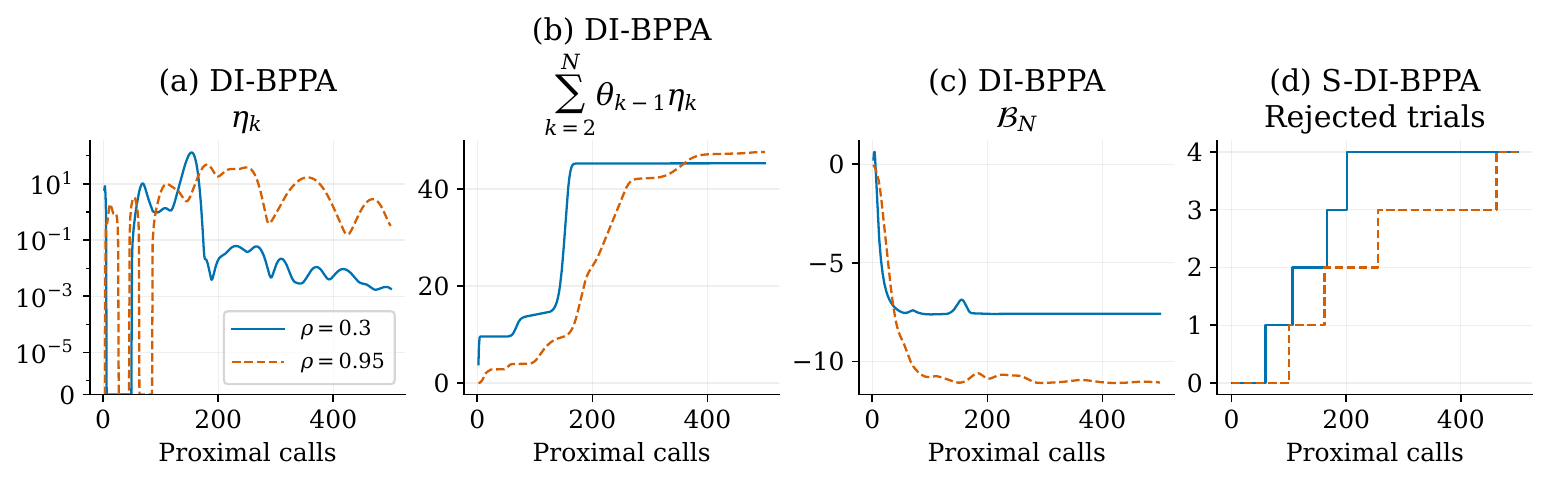}
\caption{Nonnegative least squares, $480\times320$, seed 1: geometry of unrestarted DI-BPPA in (a)--(c) and cumulative S-DI-BPPA rejections in (d), for $\rho=0.3,0.95$.}
\label{fig:nnls_geometry}
\end{figure}

Among the BPPA methods compared in Figure~\ref{fig:nnls_comparison}, DI-BPPA and S-DI-BPPA reached substantially smaller gaps at comparable proximal-call counts and update times.  S-DI-BPPA improved terminal accuracy over DI-BPPA in both configurations while retaining similar update times, with only four rejected trials in each configuration. DI-FS-BPPA also substantially outperformed standard BPPA and V-iBPPA across both configurations. 
Actually, within the $500$-call budget, DI-FS-BPPA completed $316$ outer iterations at $\rho=0.3$ and $302$ outer iterations at $\rho=0.95$. Thus, the ratio $C_N/N$ remained bounded, consistent with Theorem~\ref{thm:frac_hybrid}(2).
The connection budgets in Figure~\ref{fig:nnls_geometry} settled at negative values in both configurations. 

\section{Conclusion}
\label{sec:conclusion}
We analyzed dual-inertial Bregman proximal methods through an exact representation of the signed connection drift between dual and primal interpolation. This identity characterizes when unrestarted DI iterations maintain $\cO((1+\log N)/N^2)$ or $\cO(N^{-2})$ convergence and motivates two complementary safeguarding strategies.

Both safeguards use the same acceptance test. S-DI restarts upon rejection, guaranteeing a global $\cO(T^{-1})$ objective-gap bound under $\lambda L\le1$; under a strict stepsize bound $\lambda L<1$, local quadratic growth, and local kernel regularity at an interior minimizer, it restarts finitely often and eventually converges linearly,  while local weak sharpness yields finite restarts and finite optimal identification without $C^2$ kernel regularity.
Conversely, DI-FS preserves the quadratic weights via dual substeps. DI-FS terminates without triangle scaling and allows boundary minimizers, ensuring a global $\cO(N^{-2})$ outer-iteration rate and up to an $\cO(T^{-2})$ total oracle-call rate. Numerical experiments illustrate these theoretical bounds and demonstrate that both safeguards incur negligible computational overhead. Extending these techniques to adaptive stepsizes and inexact proximal rules remains an open direction for future work.

\begin{appendices}
\numberwithin{equation}{section}
\renewcommand{\theequation}{\thesection\arabic{equation}}
\renewcommand{\theHequation}{\thesection\arabic{equation}}

\section{Drift Representation and Directional Curvature Conditions}
\label{app:drift}
Throughout this appendix, the standing assumptions, exact DI updates, and initialization of Section~\ref{sec:analysis} apply.

\subsection{Proof of the Drift Representation}
\label{app:drift_representation}
\begin{proof}[Proof of Lemma~\ref{lem:connection_integral}]
Use the notation of the lemma and set $u'=u+pd=u_k$ and $v=u+\rho d=\nabla\omega(T_q(x_k,z_k))$. Define
\begin{equation}\label{eq:gamma_weight}
W_{p,q}(t)=
\begin{cases}
\displaystyle
\frac{t(1-p)[2p-(1+p)t]}{2p^2},
&0\le t\le p,\\[2mm]
\displaystyle
\frac{(\rho-t)^2}{2q^2},
&p\le t\le\rho.
\end{cases}
\end{equation}
Since $x_k=T_p(x_{k-1},z_k)$ and $T_q(x_k,z_k)=T_\rho(x_{k-1},z_k)$, the definition of $\Gamma_{k,q}$ and Bregman conjugacy give
\[
\Gamma_{k,q}
=D_H(s_k,u)-D_H(s_k,u')
-p^{-2}D_H(u',u)+q^{-2}D_H(v,u').
\]
Since $s_k-u'=(1-p)d$ and $u'-u=pd$, the three-point identity gives
\[
D_H(s_k,u)-D_H(s_k,u')
=p^{-1}D_H(u',u)+\frac{1-p}{p}D_H(u,u').
\]
Substitution yields
\begin{equation}\label{eq:gamma_three_divergences}
\Gamma_{k,q}
=-\frac{1-p}{p^2}D_H(u',u)
+\frac{1-p}{p}D_H(u,u')
+\frac1{q^2}D_H(v,u').
\end{equation}

The next identities require only $C^2$ regularity of $H$ near $[u,v]$. Let $m(t)=d^\top\nabla^2H(u+td)d$. Taylor's integral formula gives
\[
D_H(u+bd,u+ad)
=\int_a^b(b-t)m(t)\,dt,
\qquad a,b\in[0,\rho].
\]
Applying this formula with $(b,a)=(p,0)$, $(0,p)$, and $(\rho,p)$ in \eqref{eq:gamma_three_divergences} yields
\begin{equation}\label{eq:gamma_hessian_integral}
\begin{aligned}
\Gamma_{k,q}
=\frac{1-p}{p^2}\int_0^p[(1+p)t-p]m(t)\,dt
+\frac1{q^2}\int_p^\rho(\rho-t)m(t)\,dt
=-\int_0^\rho W_{p,q}'(t)m(t)\,dt,
\end{aligned}
\end{equation}
where $W_{p,q}'$ is taken piecewise. If $H$ is $C^3$ near $[u,v]$, then $m'(t)=\nabla^3H(u+td)[d,d,d]$. Since $W_{p,q}$ is continuous and vanishes at both endpoints, piecewise integration by parts proves \eqref{eq:gamma_integral}.

Finally, $W_{p,q}\ge0$ since $2p-(1+p)t\ge p(1-p)\ge0$ on $[0,p]$ and the second branch is nonnegative, and direct integration of \eqref{eq:gamma_weight} proves \eqref{eq:gamma_mass}.
\end{proof}

\subsection{Directional Curvature Monotonicity}
Under $C^3$ regularity, directional curvature monotonicity (DCM) is equivalent to $\eta_k=0$ for every $k\ge2$. The result below uses only the $C^2$ representation \eqref{eq:gamma_hessian_integral}.
\begin{corollary}[Directional curvature monotonicity]
\label{thm:dcm}
Assume $\lambda L\le1$. For each $k\ge2$, suppose that $H$ is $C^2$ on an open neighborhood of $[u_{k-1},s_k]$, with
\begin{equation}\label{eq:dcm}
t\longmapsto d_k^\top\nabla^2H(u_{k-1}+td_k)d_k
\end{equation}
nonincreasing on $[0,1]$. Then $\Gamma_k,r_k\le0$ for $k\ge2$, and
\[
0\le A_N\bigl(F(x_N)-F(x^*)\bigr)\le E_N\le E_1,
\qquad N\ge1.
\]
Consequently, the last-iterate objective gap is $\cO(N^{-2})$.
\end{corollary}
\begin{proof}
Fix $k\ge2$ and set $p=\theta_{k-1}$, $q=\theta_k$,
$\rho=p+q(1-p)$, and $\tau=p/(1+p)$. Define
$m(t)=d_k^\top\nabla^2H(u_{k-1}+td_k)d_k$
for $t\in[0,1]$.
The weight $W_{p,q}$ vanishes at both endpoints,
with $W_{p,q}'\ge0$ before $\tau$ and
$W_{p,q}'\le0$ after it.
Since $m$ is nonincreasing,
\eqref{eq:gamma_hessian_integral} gives
\[
\Gamma_k
=-\int_0^\rho W_{p,q}'(t)\,[m(t)-m(\tau)]\,dt
\le0.
\]

By Lemma~\ref{lem:connection_integral},
$r_k=\Gamma_{k,1}-\Gamma_{k,q}$. In
\eqref{eq:gamma_hessian_integral}, the first integral is
independent of $q$. Rescaling the second integral to $[0,1]$
therefore gives
\[
r_k=(1-p)^2\int_0^1(1-s)
\bigl[m(p+(1-p)s)-m(p+q(1-p)s)\bigr]\,ds\le0,
\]
where the sign follows from monotonicity of $m$.

Thus $\mathcal B_N\le0$ for $N\ge2$, and
\eqref{eq:drift_budget_finite} gives the stated energy and
objective-gap bounds. For $N=1$, the bounds follow from
the definition of $E_1$. The rate follows from
$A_N=\Theta(N^2)$.
\end{proof}

\subsection{Cone Conditions and Separable Kernels}
\label{app:kernel_conditions}\label{sec:problem_classes}
\begin{remark}[Sufficient geometric bounds]
\label{thm:connection_nearlog}
Under the assumptions of Theorem~\ref{prop:directional_budget}, the bound $\|\nabla^3H(y)\|\le M_3$ on the full dual segments implies $\eta_k\le M_3\|d_k\|^3$, where the tensor norm is induced by the Euclidean norm. Thus bounded $d_k$ yields the $\cO((1+\log N)/N^2)$ objective-gap rate. 
The uniform derivative bound and boundedness of $d_k$ hold if all $x_k,z_k$ lie in a compact set $K\subset C$ and $H$ is $C^3$ on an open neighborhood of $\operatorname{conv}(\nabla\omega(K))$.
\end{remark}

If the oracle vectors lie in a convex cone, the auxiliary dual update keeps the directions of the full dual segments in its negative cone.
\begin{proposition}\label{cor:dcm_cone}
Let $\mathcal K$ be a convex cone in dual coordinates. Suppose that, for every $k\ge2$, $H$ is $C^3$ on an open neighborhood of $[u_{k-1},s_k]$ and
\[
\nabla^3H(y)[d,d,d]\le0
\qquad(y\in[u_{k-1},s_k],\ d\in-\mathcal K).
\]
If the selected oracle vectors satisfy $g_k\in\mathcal K$ for all $k\ge2$, then \eqref{eq:dcm} holds and $\eta_k=0$. If also $\lambda L\le1$, then Corollary~\ref{thm:dcm} holds.
\end{proposition}
\begin{proof}
Since $x_1=z_1$, we have $s_2-u_1=-\lambda g_2/\theta_1\in-\mathcal K$. For $k\ge3$ the dual identity and auxiliary update give
\[
s_k-u_{k-1}
=(1-\theta_{k-2})(s_{k-1}-u_{k-2})
-\frac{\lambda}{\theta_{k-1}}g_k\in-\mathcal K.
\]
Induction gives $d_k=s_k-u_{k-1}\in-\mathcal K$. The cubic inequality then makes $t\mapsto d_k^\top\nabla^2H(u_{k-1}+td_k)d_k$ nonincreasing, so \eqref{eq:dcm} holds.
\end{proof}
If the cubic sign condition in Proposition~\ref{cor:dcm_cone}
holds globally, it suffices to require
$\partial F(x)\subseteq\mathcal K$ for all
$x\in C\cap\dom F$ in BPPA, or
$\nabla f(w)+\partial h(x)\subseteq\mathcal K$ for all
$w\in C$ and $x\in C\cap\dom h$ in BPG.
These domain-wide inclusions exclude isolated interior
minimizers if $F$ is finite near $x^*$ and some $v\ne0$
satisfies $\langle v,g\rangle\ge0$ for every $g\in\mathcal K$.
Indeed, for small $t>0$ and $y=x^*-tv$, they give
$g\in\partial F(y)\cap\mathcal K$ (take $w=y$ in BPG), so
$F(x^*)\ge F(y)+t\langle g,v\rangle\ge F(y)$.
Hence $y$ is also optimal. If the cubic sign condition holds
in every direction on a connected open set, oddness and
polarization give $\nabla^3H=0$, so $H$ is quadratic there.

\paragraph{Separable kernels.}
Let $\omega(x)=\sum_i\phi_i(x_i)$ satisfy Assumption~\ref{ass:omega}, where $I_i=\operatorname{int}(\dom\phi_i)$, $\phi_i\in C^3(I_i)$, $\phi_i''>0$, $\phi_i'''\le0$, and $\phi_i'(I_i)=\R$. Then $H(u)=\sum_i h_i(u_i)$ with $h_i=\phi_i^*$, and differentiating $h_i'=(\phi_i')^{-1}$ gives
\[
h_i'''(u)
=-\frac{\phi_i'''(x)}{[\phi_i''(x)]^3}\ge0,
\qquad x=h_i'(u).
\]
Thus $\nabla^3H(u)[d,d,d]\le0$ for $d\in-\R_+^n$, and $g_k\in\R_+^n$ for all $k\ge2$ implies DCM by Proposition~\ref{cor:dcm_cone}.
Examples include the orthant entropy kernel $\omega(x)=\sum_i(x_i\log x_i-x_i)$ and the quadratically regularized Burg kernel $\omega(x)=\gamma_2\|x\|^2/2-\gamma_1\sum_i\log x_i$, with $\gamma_1,\gamma_2>0$~\citep{bolte2003barrier}.

\section{Local Stability and Linear Convergence}
\label{app:local_analysis}
We consider exact DI updates with a fixed stepsize $\lambda>0$ and weights satisfying~\eqref{eq:param}. The local growth condition is Assumption~\ref{ass:local_qg}.

\subsection{Local Stability and Contraction of DI}
\label{ssec:bppa_refinement}
Under Assumption~\ref{ass:local_qg}, $x^*$ is a strong local minimizer
\citep{drusvyatskiy2013tilt}.
For any DI step with $x_k,x_{k+1}\in B_r(x^*)$, the definitions of $U_k$ and $V_k$ give
\[
\frac{U_k}{A_k}-V_k
=F(x_k)-F(x^*)
-\langle g_{k+1},x_k-x^*\rangle.
\]
Since $V_k\ge0$, applying \eqref{eq:local_qg} and Young's inequality to this identity and the definition of $V_k$ yields
\begin{equation}\label{eq:local_oracle_gap}
\begin{aligned}
F(x_k)-F(x^*)
&\le\frac{2U_k}{A_k}
+\frac{2}{m_F}\|g_{k+1}\|^2,\\
F(x_{k+1})-F(x^*)
&\le\frac{2}{m_F}\|g_{k+1}\|^2
+2D_f(x_{k+1},w_k).
\end{aligned}
\end{equation}

We next establish local stability of the DI methods. The kernel assumptions are satisfied by the Euclidean, orthant entropy, and quadratically regularized Burg kernels, which are $C^2$ on $C$ with positive definite Hessians.
\begin{theorem}[Local stability of the DI methods]
\label{thm:local_bppa_stability}
Under Assumption~\ref{ass:local_qg}, suppose that $\lambda L<1$ and that $\omega$ is $C^2$ near $x^*$ with $\nabla^2\omega(x^*)\succ0$. There exist $\varepsilon_*>0$ and $c>0$ such that any exact DI method with weights satisfying \eqref{eq:param} and initial state $x_1\in\dom F\cap C$, $z_1\in C$, satisfying
\[
E_1=A_1(F(x_1)-F(x^*))+D_\omega(x^*,z_1)
\le\varepsilon_*,
\]
keeps all $x_k,z_k,w_k$ in a fixed compact subset of $C$ and satisfies
\begin{equation}\label{eq:local_bppa_step}
E_{k+1}+cD_\omega(z_{k+1},z_k)+\tfrac12S_k\le E_k,
\qquad k\ge1.
\end{equation}
The constants can be chosen independently of $A_1$ when $A_1\ge\lambda$. For DI-BPPA, every fixed $\lambda>0$ is allowed, and one may take $c=1/2$.
\end{theorem}

\begin{proof}
Let $\sigma=1-\lambda L>0$.  Consider  $\mu I\preceq\nabla^2\omega\preceq MI$ on a convex primal neighborhood contained in $B_r(x^*)$ and the corresponding dual bounds $M^{-1}I\preceq\nabla^2H\preceq\mu^{-1}I$ near $u^*=\nabla\omega(x^*)$.
These bounds remain fixed as the neighborhoods shrink. The convexity of $L\omega-f$ and $(M/2)\|\cdot\|^2-\omega$ implies that $\nabla f$ is locally $L_f$-Lipschitz with $L_f=LM$. For simplicity, throughout the proof, $K$ denotes positive constants depending only on these fixed bounds and the problem parameters.

First consider a state with sufficiently small energy $E_k$. Choose $R\in(0,r)$ sufficiently small. 
For $y\in\dom F$ with $\|y-x^*\|\ge R$, convexity along $[x^*,y]$ and \eqref{eq:local_qg} imply
\[
F(y)-F(x^*)
\ge\frac{\|y-x^*\|}{R}
\left[
F\!\left(x^*+\frac{R(y-x^*)}{\|y-x^*\|}\right)-F(x^*)
\right]\ge\frac{m_FR}{2}\|y-x^*\|.
\]
This bound and Lemma~\ref{lem:bregman_sublevel_compact} place $x_k$ and $s_k$ near $x^*$ and $u^*$, respectively, with
\begin{equation}\label{eq:local_current_bounds}
\|x_k-x^*\|\le\sqrt{\frac{2E_k}{m_FA_k}},
\qquad
\|s_k-u^*\|\le\sqrt{2ME_k}.
\end{equation}
Suppress the index $k$, writing $A=A_k$, $\theta=\theta_k$,
$x=x_k$, $p=x_{k+1}$, $u=u_k$, $u^+=u_{k+1}$,
$s=s_k$, $s^+=s_{k+1}$, and $g=g_{k+1}$.
Let $v=(1-\theta)u+\theta s$ and $w=\nabla H(v)$.
By Lemma~\ref{lem:prox_interior} with $\varphi=h$, the oracle output $p$ depends continuously on $t=v-\lambda\nabla f(w)$, and hence on $v$.
At $v=u^*$, the optimality condition gives $p=x^*$.
So $p$ lies in the primal neighborhood whenever $v$ is sufficiently close to $u^*$.
The optimality conditions and monotonicity of $\partial h$ give
\[
\mu\|p-x^*\|^2
\le
\left\langle
v-u^*-\lambda\bigl(\nabla f(w)-\nabla f(x^*)\bigr),
p-x^*
\right\rangle.
\]
Since $\|w-x^*\|\le\mu^{-1}\|v-u^*\|$ and $\nabla f$ is $L_f$-Lipschitz there, Cauchy--Schwarz and the upper Hessian bound yield
\[
\|p-x^*\|+\|u^+-u^*\|\le K\|v-u^*\|.
\]
Reconstruction, together with
$(1-\theta)/(\theta\sqrt A)=\sqrt{(1-\theta)/\lambda}$, yields
\[
\|s^+-u^*\|
\le K\left(
\frac{1-\theta}{\theta}\sqrt{\frac{E_k}{A}}
+\sqrt{E_k}\right)
\le K\sqrt{E_k}.
\]
Thus sufficiently small current energy places all trial points and relevant segments in any prescribed sufficiently small neighborhoods.

Let $b=u^+-u$, $\Delta s=s^+-s$, $D=D_\omega(z_{k+1},z_k)$, and $\widehat Q=\widehat Q_k$. Since $v-u^+=-\theta\Delta s$, the dual updates and local Hessian bounds imply
\[
\begin{aligned}
\widehat Q
&\ge\sigma\theta^{-2}D_H(v,u^+)
\ge\frac{\sigma}{2M}\|\Delta s\|^2,\\
\|g\|^2
&=\frac{\theta^2}{\lambda^2}\|\Delta s\|^2,
\qquad
D_f(p,w)\le\frac{L\theta^2}{2\mu}\|\Delta s\|^2.
\end{aligned}
\]
By \eqref{eq:local_qg} and \eqref{eq:local_oracle_gap},
\[
\|p-x\|^2
\le2\|p-x^*\|^2+2\|x-x^*\|^2
\le\frac{8U_k}{m_FA}
+\frac{16}{m_F^2}\|g\|^2
+\frac{8}{m_F}D_f(p,w).
\]
Using $A\theta^2\le\lambda$, we obtain
\[
A\|p-x\|^2
\le\frac{8}{m_F}U_k
+\left(
\frac{16}{m_F^2\lambda}
+\frac{4\lambda L}{m_F\mu}
\right)\|\Delta s\|^2.
\]
Together with $\|b\|\le M\|p-x\|$ and the preceding lower bound for $\widehat Q$, we obtain
\begin{equation}\label{eq:local_displacement_control}
\|\Delta s\|^2+A\|b\|^2
\le K(\widehat Q+U_k).
\end{equation}

It remains to bound the connection increment. For fixed $u,\theta$, define
\[
\widetilde{\mathcal R}_\theta(t)
:=D_H(t,u)-\theta^{-2}D_H(u+\theta(t-u),u).
\]
Let $P=\nabla^2H(u^*)$ and shrink the dual neighborhood so that $\|\nabla^2H-P\|\le\delta$. With $q=t-u$, differentiation gives
\[
\nabla\widetilde{\mathcal R}_\theta(t)
=\int_\theta^1
[\nabla^2H(u+\tau q)-P]q\,d\tau
-\frac{1-\theta}{\theta}
\int_0^\theta[\nabla^2H(u+\tau q)-P]q\,d\tau.
\]
Hence $\|\nabla\widetilde{\mathcal R}_\theta(t)\|
\le2\delta(1-\theta)\|t-u\|$. 
In these dual coordinates, the connection increment is
\[
J_k=\widetilde{\mathcal R}_\theta(s^+)
    -\widetilde{\mathcal R}_\theta(s).
\]
Since $s^+-u=b/\theta$, integration along $[s,s^+]$ gives
\[
|J_k|
\le2\delta\frac{1-\theta}{\theta}\|\Delta s\|\,\|b\|
+\delta(1-\theta)\|\Delta s\|^2
\le K\delta(\widehat Q+U_k),
\]
where the last step uses \eqref{eq:local_displacement_control} and
$(1-\theta)/\theta\le\sqrt{A/\lambda}$.

With $K$ as in the preceding estimate, choose $\delta>0$
so that $K\delta\le1/4$. By continuity and positive definiteness, the dual neighborhood may also be chosen so that its Hessians are pairwise comparable within a factor $3/2$. Taylor's integral formula then gives
\[
\widehat Q
\ge\sigma\theta^{-2}D_H(v,u^+)
\ge\frac{2\sigma}{3}D.
\]
Since $J_k=\Gamma_{k+1}+r_{k+1}-r_k$, the exact balance \eqref{eq:exact_connection_recursion} becomes
\[
E_{k+1}-E_k=-\widehat Q-S_k+J_k.
\]
Together with $U_k\le S_k$, this proves
\[
E_{k+1}-E_k
\le-\tfrac34\widehat Q-\tfrac34S_k
\le-\frac{\sigma}{2}D-\frac12S_k.
\]
The same estimates give the margin bound
\begin{equation}\label{eq:local_margin_control}
\widehat Q_k+U_k-J_k
\ge\frac{\sigma}{2}D_\omega(z_{k+1},z_k)+\frac12U_k.
\end{equation}

Choose $\varepsilon_*$ small enough for the preceding localization. Induction proves \eqref{eq:local_bppa_step} with $c=\sigma/2$ and compact containment. 
All estimates are uniform for $A_1\ge\lambda$, since $A_k\ge A_1$ and the weighted reconstruction bounds do not require a lower bound on $1-\theta_k$.
\end{proof}
{ 

\begin{corollary}[Local linear convergence of the DI methods]
\label{cor:local_bppa_linear}
Under the assumptions and initialization of Theorem~\ref{thm:local_bppa_stability}, there exists $\rho\in(0,1)$ such that, for every $k\ge1$,
\begin{equation}\label{eq:local_bppa_linear}
E_k\le E_1\rho^{k-1},\qquad
0\le F(x_k)-F(x^*)\le\frac{E_1}{A_k}\rho^{k-1}.
\end{equation}
In particular, $\|x_k-x^*\|=\cO(\rho^{k/2}/k)$ and $\|z_k-x^*\|=\cO(\rho^{k/2})$. 
In addition, if $F$ is locally $m_F$-strongly convex, then for DI-BPPA one may take $\rho=\left(1+\frac{\mu m_F\lambda}{2M^2}\right)^{-1}$, where $\mu I\preceq\nabla^2\omega\preceq MI$.
\end{corollary}
\begin{proof}
We use the local Hessian bounds from the preceding proof,
with $K$ denoting positive constants independent of $k$,
and put $\sigma=1-\lambda L>0$.
The second inequality in \eqref{eq:local_oracle_gap},
relative smoothness, and the local dual Hessian bounds give
\[
F(x_{k+1})-F(x^*)
\le
\left(\frac{2}{m_F}+\frac{\lambda^2L}{\mu}\right)
\|g_{k+1}\|^2.
\]
Moreover,
\begin{equation*}
\widehat Q_k
\ge
\frac{\sigma\lambda^2}{2M\theta_k^2}
\|g_{k+1}\|^2, \qquad
\|g_{k+1}\|^2
=\frac{\theta_k^2}{\lambda^2}
\|s_{k+1}-s_k\|^2
\le
\frac{2M\theta_k^2}{\lambda^2}
D_\omega(z_{k+1},z_k).
\end{equation*}
Using $A_{k+1}\theta_k^2=\lambda$, we obtain
\begin{equation}\label{eq:local_gap_dissipation}
\begin{aligned}
A_{k+1}\bigl(F(x_{k+1})-F(x^*)\bigr)
&\le K\widehat Q_k,\\
A_{k+1}\bigl(F(x_{k+1})-F(x^*)\bigr)
&\le K D_\omega(z_{k+1},z_k).
\end{aligned}
\end{equation}
The first estimate only requires $x_{k+1}$ and $w_k$
to lie in the local neighborhoods.
The local Hessian bounds and $D_f\ge0$ give
$\widehat Q_k\le(M/\mu)D_\omega(z_{k+1},z_k)$.
Hence \eqref{eq:local_displacement_control} implies
\[
A_k\|u_{k+1}-u_k\|^2
\le K\bigl[D_\omega(z_{k+1},z_k)+U_k\bigr].
\]
The  bound \eqref{eq:local_qg} and the upper  Hessian bound also imply
$\|u_{k+1}-u^*\|^2
\le K[F(x_{k+1})-F(x^*)]$.
The weighted reconstruction identity yields
\[
\begin{aligned}
\lambda\|s_{k+1}-u^*\|^2
&\le A_k\|u_{k+1}-u_k\|^2
+a_{k+1}\|u_{k+1}-u^*\|^2\\
&\le K\bigl[D_\omega(z_{k+1},z_k)+S_k\bigr],
\end{aligned}
\]
where the last inequality uses the second estimate in
\eqref{eq:local_gap_dissipation},
$a_{k+1}\le A_{k+1}$, and $U_k\le S_k$.
Since
$D_\omega(x^*,z_{k+1})
\le\|s_{k+1}-u^*\|^2/(2\mu)$,
the second estimate in \eqref{eq:local_gap_dissipation}
now gives
\[
E_{k+1}
\le K\bigl[D_\omega(z_{k+1},z_k)+S_k\bigr]
\le\frac{K}{c}(E_k-E_{k+1}),
\]
where the last inequality follows from
\eqref{eq:local_bppa_step} with
$c=(1-\lambda L)/2\le1/2$.
Thus \eqref{eq:local_bppa_linear} holds with
$\rho=(1+c/K)^{-1}$.

Suppose additionally that $F$ is locally $m_F$-strongly convex. For DI-BPPA, $D_f=0$ and $g_{k+1}\in\partial F(x_{k+1})$. The strong-convexity inequalities, the auxiliary update, and the local kernel Hessian bounds give
\[
\begin{aligned}
A_{k+1}(F(x_{k+1})-F(x^*))
&\le \frac{M^2}{\mu m_F\lambda}
D_\omega(z_{k+1},z_k),\\
D_\omega(x^*,z_{k+1})
&\le \frac{M^2}{\mu m_F\lambda}S_k.
\end{aligned}
\]
For DI-BPPA, $L=0$ and hence $c=1/2$, so
\[
E_{k+1}
\le \frac{M^2}{\mu m_F\lambda}
\bigl[D_\omega(z_{k+1},z_k)+S_k\bigr]
\le \frac{2M^2}{\mu m_F\lambda}(E_k-E_{k+1}),
\]
which yields the stated BPPA contraction factor.
For $A_1\ge\lambda$, the uniform local estimates make
$\rho$ independent of $A_1$.

Finally, combining \eqref{eq:local_current_bounds} with $E_k\le E_1\rho^{k-1}$, $A_k=\Theta(k^2)$, and $\|z_k-x^*\|=\|\nabla H(s_k)-\nabla H(u^*)\|\le\mu^{-1}\|s_k-u^*\|$ yields $\|x_k-x^*\|=\cO(\rho^{k/2}/k)$ and $\|z_k-x^*\|=\cO(\rho^{k/2})$.
\end{proof}

\begin{remark}[Strong convexity and quadratic growth]
At an interior minimizer, local $m_F$-strong convexity implies Assumption~\ref{ass:local_qg}. Likewise, local relative strong convexity with modulus $m_{\mathrm{rel}}>0$ gives
$
F(x)-F(x^*)
\ge m_{\mathrm{rel}}D_\omega(x,x^*)
\ge\frac{m_{\mathrm{rel}}\mu}{2}\|x-x^*\|^2,
$
where $\nabla^2\omega\succeq\mu I$ locally. Thus Assumption~\ref{ass:local_qg} holds with $m_F=m_{\mathrm{rel}}\mu$.
\end{remark}

\subsection{Proof of Eventual Linear Convergence}
\label{app:eventual_linear_proof}
\begin{proof}[Proof of Corollary~\ref{cor:bppa_eventual_linear}]

It suffices to prove $E_j\to0$. Once $E_j\le\varepsilon_*$,
the uniform estimates of
Theorem~\ref{thm:local_bppa_stability} and
Corollary~\ref{cor:local_bppa_linear} apply with the current
weight $A_j\ge\lambda$.

All trials in the final epoch are accepted, so
\eqref{eq:adaptive_edge_energy} gives $0\le E_j\le E_1$.
Set $\phi_j:=F(x_j)-F(x^*)$ and $u^*:=\nabla\omega(x^*)$.
Since $A_j=\Theta(j^2)$, we have $\phi_j=\cO(j^{-2})$.
The radial lower bound in the proof of
Theorem~\ref{thm:local_bppa_stability} then gives $x_j\to x^*$.
Moreover,
$D_H(s_j,u^*)=D_\omega(x^*,z_j)\le E_1$,
so Lemma~\ref{lem:bregman_sublevel_compact} bounds $s_j$.
Consequently,
$v_j:=\nabla\omega(w_j)=(1-\theta_j)u_j+\theta_js_j\to u^*$,
and hence $w_j\to x^*$.

Consider $\mu I\preceq\nabla^2\omega\preceq MI$ from Theorem~\ref{thm:local_bppa_stability}, and set
\[
P:=\nabla^2H(u^*)\succ0,\qquad
\|b\|_P^2:=b^\top Pb,\qquad
b_j:=u_j-u_{j-1},\qquad
\beta_j:=\frac{1-\theta_{j-1}}{\theta_{j-1}}.
\]
Reconstruction gives $s_j=u_j+\beta_jb_j$ and $v_j-u_j=\theta_j\beta_jb_j$. The three-point identity yields
\[
\Phi_j
=A_j\phi_j+D_\omega(x^*,x_j)
+\theta_j^{-2}D_H(v_j,u_j)
+\beta_j\langle b_j,x_j-x^*\rangle.
\]
For sufficiently large $j$, the local Hessian bounds, \eqref{eq:local_qg}, and Young's inequality give
\[
\begin{aligned}
\Phi_j
&\ge\left(A_j-\frac{2M}{m_F}\right)\phi_j
+\frac{\beta_j^2}{4M}\|b_j\|^2
\ge\frac12A_j\phi_j+\frac{\beta_j^2}{4M}\|b_j\|^2,\\
\Phi_j
&\le K\bigl(A_j\phi_j+\beta_j^2\|b_j\|^2\bigr),
\end{aligned}
\]
where $K$ is independent of $j$ and the second lower bound uses $A_j\to\infty$.

Since $x_{j-1},x_j,w_{j-1}\to x^*$, the derivation of
\eqref{eq:local_displacement_control} applies for all
sufficiently large $j$. Equivalence of the Euclidean
and $P$-norms then gives
\[
A_{j-1}\|b_j\|_P^2\le K(\widehat Q_{j-1}+U_{j-1}).
\]
Since $x_j,w_{j-1}\to x^*$, the first estimate in
\eqref{eq:local_gap_dissipation} applies at iteration $j-1$
for all sufficiently large $j$, giving
$A_j\phi_j\le K\widehat Q_{j-1}$.
Together with $\beta_j^2\le A_{j-1}/\lambda$, these estimates allow us to fix $C_0\ge1$ such that
\[
\Phi_j\le C_0(\widehat Q_{j-1}+U_{j-1}),
\qquad
\lambda^{-1}A_{j-1}\|b_j\|_P^2
\le C_0(\widehat Q_{j-1}+U_{j-1})
\]
for all sufficiently large $j$.

Choose $0<\epsilon\le1/2$ with $\epsilon C_0\le1/2$. By continuity, a smaller convex dual neighborhood of $u^*$ satisfies
\[
(1-\epsilon)P\preceq\nabla^2H\preceq(1+\epsilon)P.
\]
The constant $C_0$ remains fixed, and $u_{j-1},u_j,v_j$ eventually lie in this neighborhood. Apply the resulting quadratic bounds to the three divergences in \eqref{eq:gamma_three_divergences}. With $p=\theta_{j-1}$ and $v_j-u_j=\theta_j\beta_jb_j$, this gives
\[
\begin{aligned}
\Gamma_j
\le\frac{1-p}{2p^2}
\bigl[-(1-\epsilon)+p(1+\epsilon)+(1-p)(1+\epsilon)\bigr]
\|b_j\|_P^2
=\frac{\epsilon A_{j-1}}{\lambda}\|b_j\|_P^2
\le\frac12(\widehat Q_{j-1}+U_{j-1}).
\end{aligned}
\]
Thus \eqref{eq:exact_connection_recursion} and $V_{j-1}\ge0$ imply
\[
\Phi_j+\frac12(\widehat Q_{j-1}+U_{j-1})
\le\Phi_{j-1},
\qquad
\left(1+\frac1{2C_0}\right)\Phi_j\le\Phi_{j-1}.
\]
The lower bound on $\Phi_j$ now gives $A_j\phi_j\to0$ and $\beta_jb_j\to0$. Since $s_j=u_j+\beta_jb_j$ and $u_j\to u^*$, we have $s_j\to u^*$. Continuity of $H$ therefore yields
\[
E_j=A_j\phi_j+D_H(s_j,u^*)\longrightarrow0.
\]
Choose $j_0$ with $E_{j_0}\le\varepsilon_*$.
The local contraction estimate and
\eqref{eq:weight_lower_bound} prove
\eqref{eq:bppa_eventual_linear}.
\end{proof}

\begin{remark}[Local regularity without global quadratic triangle scaling]
\label{app:burg_scaling_separation}
The quadratically regularized Burg kernel from
Appendix~\ref{app:kernel_conditions} does not satisfy global
quadratic triangle scaling. In one dimension, its divergence is
$D_\omega(r,t)=\gamma_1(r/t-1-\log(r/t))+\gamma_2(r-t)^2/2$.
In the proposed scaling inequality
\[
D_\omega((1-\theta)b+\theta x,(1-\theta)b+\theta y)
\le G\theta^2D_\omega(x,y),
\]
take $b=\theta^2$, $x=2$, and $y=1$. As $\theta\downarrow0$,
the ratio of the two arguments on the left tends to $2$
and their difference tends to zero. Thus the left-hand side
tends to $\gamma_1(1-\log2)>0$, whereas the right-hand side
tends to zero for every finite $G$. Fixing the remaining
coordinates gives the same counterexample in higher dimensions.

The kernel is nevertheless $C^\infty$ on $C$, with
$\nabla^2\omega(x)=\gamma_2I+\gamma_1\operatorname{diag}(x_i^{-2})
\succ0$. Thus it meets the kernel regularity requirement of
Theorem~\ref{thm:safeguard_finite_restart}.
\end{remark}
}

\section{Proof of Example~\ref{rem:kl_quadratic}}
\label{app:kl_quadratic}

Apply Lemma~\ref{lem:prox_interior} to the indicator function of the affine set $\{x:Bx=b\}$ at $x_1$. The projection $x^\dagger$ exists uniquely in $C$, and its optimality condition gives $u^\dagger-u_1\in\operatorname{Range}(B)$. Choose $a\in C$ with $Ba=b$. Then
\[
F(x)-F^*=\tfrac12(x-a)^\top B(x-a).
\]
Using the dual coordinates of Section~\ref{ssec:connection_drift}, define
$G(u):=H(u)-\langle a,u\rangle$,
$\mathcal S:=u_1+\operatorname{Range}(B)$, and
$v_k:=(1-\theta_k)u_k+\theta_ks_k$.
Since $Bx^\dagger=b$ and $u^\dagger\in\mathcal S$, convexity makes $u^\dagger$ a minimizer of $G$ on $\mathcal S$. By orthogonality and \eqref{eq:bregman_dual}, we have,  for $u\in\mathcal S$,
\[
G(u)-G(u^\dagger)
=D_H(u,u^\dagger)
=D_\omega(x^\dagger,\nabla H(u)).
\]
Proposition~\ref{prop:common_dual_update} gives
\[
u_{k+1}=v_k-\lambda B(\nabla H(u_{k+1})-a),
\qquad
s_{k+1}=s_k-\frac{\lambda}{\theta_k}B(\nabla H(u_{k+1})-a).
\]
Since $s_1=u_1$, all dual iterates lie in $\mathcal S$. The first recursion is equivalent to $\nabla G(u_{k+1})+\lambda^{-1}B^\dagger(u_{k+1}-v_k)\in\ker B$. As $B^\dagger$ is positive definite on $\operatorname{Range}(B)$,
\begin{equation}\label{eq:kl_quadratic_metric_prox}
u_{k+1}
=\arg\min_{u\in\mathcal S}
\left\{G(u)+\frac{1}{2\lambda}\|u-v_k\|_{B^\dagger}^2\right\}.
\end{equation}
Together with affine interpolation and reconstruction, this is
a fixed-metric accelerated proximal point iteration
\citep[Section~3]{monteiro2013accelerated}.

Since $BB^\dagger B=B$ and $A_{k+1}\theta_k^2=\lambda$,
\[
\frac12\|s_{k+1}-s_k\|_{B^\dagger}^2
=\frac{\lambda^2}{2\theta_k^2}
(x_{k+1}-a)^\top B(x_{k+1}-a)
=\lambda A_{k+1}\bigl(F(x_{k+1})-F^*\bigr).
\]
Apply \eqref{eq:lyap_descent_bppa} to $G$ on $\mathcal S$
with the quadratic kernel
$\tfrac12\|\cdot-u_1\|_{B^\dagger}^2$.
The interpolation defect vanishes, so summation gives, for $N\ge2$,
\begin{equation*}\label{eq:kl_quadratic_energy}
A_ND_\omega(x^\dagger,x_N)
+\tfrac12\|s_N-u^\dagger\|_{B^\dagger}^2
+\lambda\sum_{j=2}^N A_j\bigl(F(x_j)-F^*\bigr)
\le C_{\rm init}.
\end{equation*}
Nonnegativity proves both stated bounds.

\end{appendices}

\section*{Statements and Declarations}

\noindent\textbf{Funding}  
The author received no external funding for this work.

\noindent\textbf{Competing interests}
The author declares no competing interests.

\noindent\textbf{Data and code availability}
Numerical instances were generated as described in Section~\ref{sec:experiments_numeric}.

\noindent\textbf{Use of generative AI}
The author used OpenAI's ChatGPT and Codex to assist with
manuscript drafting and language editing, literature searches,
mathematical derivations and proof development, and the
implementation of numerical experiments.
The author takes full responsibility for the content.

\bibliography{references}

@article{chen2026accmd,
  author = {Chen, Long and Luo, Hao and Wei, Jingrong and Xu, Zeyi and Yao, Yuan},
  title = {Accelerated Mirror Descent Method through Variable and Operator Splitting},
  journal = {arXiv preprint arXiv:2601.19038},
  year = {2026}
}

@article{gutman2018unified,
  title   = {A unified framework for {Bregman} proximal methods: subgradient, gradient, and accelerated gradient schemes},
  author  = {Gutman, David H. and Pe{\~n}a, Javier F.},
  journal = {arXiv preprint arXiv:1812.10198},
  year    = {2018}
}

@article{yan2020bregman,
  author  = {Yan, Shen and He, Niao},
  title   = {{Bregman} Augmented {Lagrangian} and Its Acceleration},
  journal = {arXiv preprint arXiv:2002.06315},
  year    = {2020}
}

@article{liu2022dual,
  author  = {Liu, Jin-Zan and Liu, Xin-Wei},
  title   = {A dual {Bregman} proximal gradient method for relatively-strongly convex optimization},
  journal = {Numer. Algebra Control Optim.},
  volume  = {12},
  number  = {4},
  pages   = {679--692},
  year    = {2022},
  doi     = {10.3934/naco.2021028}
}

@article{azizian2024rate,
  author  = {Azizian, Wa{\"i}ss and Iutzeler, Franck and Malick, J{\'e}r{\^o}me and Mertikopoulos, Panayotis},
  title   = {The Rate of Convergence of {Bregman} Proximal Methods: Local Geometry Versus Regularity Versus Sharpness},
  journal = {SIAM J. Optim.},
  volume  = {34},
  number  = {3},
  pages   = {2440--2471},
  year    = {2024}
}

@misc{allen2016optimal,
  title         = {Linear Coupling: An Ultimate Unification of Gradient and Mirror Descent},
  author        = {Zeyuan Allen-Zhu and Lorenzo Orecchia},
  year          = {2016},
  eprint        = {1407.1537},
  archiveprefix = {arXiv},
  primaryclass  = {cs.DS}
}

@article{bauschke2017,
  author  = {Bauschke, Heinz H. and Bolte, J{\'e}r{\^o}me and Teboulle, Marc},
  title   = {A Descent Lemma Beyond {Lipschitz} Gradient Continuity: First-Order Methods Revisited and Applications},
  journal = {Math. Oper. Res.},
  volume  = {42},
  number  = {2},
  pages   = {330--348},
  year    = {2017}
}

@article{beck2003mirror,
  author  = {Beck, Amir and Teboulle, Marc},
  title   = {Mirror Descent and Nonlinear Projected Subgradient Methods for Convex Optimization},
  journal = {Oper. Res. Lett.},
  volume  = {31},
  number  = {3},
  pages   = {167--175},
  year    = {2003}
}

@article{beck2009fista,
  author  = {Beck, Amir and Teboulle, Marc},
  title   = {A Fast Iterative Shrinkage-Thresholding Algorithm for Linear Inverse Problems},
  journal = {SIAM J. Imaging Sci.},
  volume  = {2},
  number  = {1},
  pages   = {183--202},
  year    = {2009}
}

@article{beck2014fast,
  author  = {Beck, Amir and Teboulle, Marc},
  title   = {A Fast Dual Proximal Gradient Algorithm for Convex Minimization and Applications},
  journal = {Oper. Res. Lett.},
  volume  = {42},
  number  = {1},
  pages   = {1--6},
  year    = {2014}
}

@article{bolte2003barrier,
  author  = {Bolte, J{\'e}r{\^o}me and Teboulle, Marc},
  title   = {Barrier Operators and Associated Gradient-Like Dynamical Systems for Constrained Minimization Problems},
  journal = {SIAM J. Control Optim.},
  volume  = {42},
  number  = {4},
  pages   = {1266--1292},
  year    = {2003},
  doi     = {10.1137/S0363012902410861}
}

@article{bregman1967,
  author  = {Bregman, Lev M.},
  title   = {The Relaxation Method of Finding the Common Point of Convex Sets and Its Application to the Solution of Problems in Convex Programming},
  journal = {U.S.S.R. Comput. Math. Math. Phys.},
  volume  = {7},
  number  = {3},
  pages   = {200--217},
  year    = {1967}
}

@article{censor1992,
  author  = {Censor, Yair and Zenios, Stavros},
  title   = {Proximal Minimization Algorithm with D-Functions},
  journal = {J. Optim. Theory Appl.},
  volume  = {73},
  number  = {3},
  pages   = {451--464},
  year    = {1992}
}

@article{chen1993,
  author  = {Chen, Gong and Teboulle, Marc},
  title   = {Convergence Analysis of a Proximal-Like Minimization Algorithm Using {Bregman} Functions},
  journal = {SIAM J. Optim.},
  volume  = {3},
  number  = {3},
  pages   = {538--543},
  year    = {1993}
}

@article{dragomir2021,
  author  = {Dragomir, Radu-Alexandru and Taylor, Adrien B. and d'Aspremont, Alexandre and Bolte, J{\'e}r{\^o}me},
  title   = {Optimal Complexity and Certification of {Bregman} First-Order Methods},
  journal = {Math. Program.},
  volume  = {194},
  number  = {1},
  pages   = {41--83},
  year    = {2022}
}

@inproceedings{defazio2019curved,
  author    = {Defazio, Aaron},
  title     = {On the Curved Geometry of Accelerated Optimization},
  booktitle = {Advances in Neural Information Processing Systems},
  volume    = {32},
  year      = {2019}
}

@article{eckstein1993,
  author  = {Eckstein, Jonathan},
  title   = {Nonlinear Proximal Point Algorithms Using {Bregman} Functions, with Applications to Convex Programming},
  journal = {Math. Oper. Res.},
  volume  = {18},
  number  = {1},
  pages   = {202--226},
  year    = {1993}
}

@article{lu2018relatively,
  author  = {Lu, Haihao and Freund, Robert M. and Nesterov, Yurii},
  title   = {Relatively Smooth Convex Optimization by First-Order Methods, and Applications},
  journal = {SIAM J. Optim.},
  volume  = {28},
  number  = {1},
  pages   = {333--354},
  year    = {2018}
}

@article{yang2022bregman,
  author  = {Yang, Lei and Toh, Kim-Chuan},
  title   = {Bregman Proximal Point Algorithm Revisited: A New Inexact Version and Its Inertial Variant},
  journal = {SIAM J. Optim.},
  volume  = {32},
  number  = {3},
  pages   = {1523--1554},
  year    = {2022}
}

@article{hanzely2021accelerated,
  author  = {Hanzely, Filip and Richt{\'a}rik, Peter and Xiao, Lin},
  title   = {Accelerated {Bregman} Proximal Gradient Methods for Relatively Smooth Convex Optimization},
  journal = {Comput. Optim. Appl.},
  volume  = {79},
  number  = {2},
  pages   = {405--440},
  year    = {2021}
}

@article{yang2026inexact,
  title     = {Inexact {Bregman} proximal gradient method and its inertial variant with absolute and partial relative stopping criteria},
  author    = {Yang, Lei and Toh, Kim-Chuan},
  journal   = {Math. Oper. Res.},
  volume    = {51},
  number    = {1},
  pages     = {686--714},
  year      = {2026},
  publisher = {INFORMS}
}

@article{nesterov2007dual,
  title     = {Dual extrapolation and its applications to solving variational inequalities and related problems},
  author    = {Nesterov, Yurii},
  journal   = {Math. Program.},
  volume    = {109},
  number    = {2},
  pages     = {319--344},
  year      = {2007},
  publisher = {Springer}
}

@book{rockafellar1970convex,
  title     = {Convex Analysis},
  author    = {Rockafellar, R. Tyrrell},
  year      = {1970},
  publisher = {Princeton University Press},
  address   = {Princeton, NJ}
}

@article{maddison2021dual,
  title     = {Dual space preconditioning for gradient descent},
  author    = {Maddison, Chris J and Paulin, Daniel and Teh, Yee Whye and Doucet, Arnaud},
  journal   = {SIAM J. Optim.},
  volume    = {31},
  number    = {1},
  pages     = {991--1016},
  year      = {2021},
  publisher = {SIAM}
}

@article{kim2023mirror,
  title   = {Mirror duality in convex optimization},
  author  = {Kim, Jaeyeon and Park, Chanwoo and Ozdaglar, Asuman and Diakonikolas, Jelena and Ryu, Ernest K},
  journal = {arXiv preprint arXiv:2311.17296},
  year    = {2023}
}

@article{monteiro2013accelerated,
  title     = {An accelerated hybrid proximal extragradient method for convex optimization and its implications to second-order methods},
  author    = {Monteiro, Renato DC and Svaiter, Benar Fux},
  journal   = {SIAM J. Optim.},
  volume    = {23},
  number    = {2},
  pages     = {1092--1125},
  year      = {2013},
  publisher = {SIAM}
}

@article{xu2026adaptive,
  author = {Xu, Zeyi and Chen, Long},
  title = {Adaptive Accelerated Mirror Descent in Primal and Dual Spaces},
  journal = {arXiv preprint arXiv:2606.02787},
  year = {2026}
}

@article{bauschke1997legendre,
  author  = {Bauschke, Heinz H. and Borwein, Jonathan M.},
  title   = {Legendre Functions and the Method of Random {Bregman} Projections},
  journal = {J. Convex Anal.},
  volume  = {4},
  number  = {1},
  pages   = {27--67},
  year    = {1997}
}

@misc{chen2026bregman,
  title         = {A {Bregman} inertial iteratively regularized extragradient method for bilevel variational inequality problems},
  author        = {Kangming Chen and Atsushi Hori and Ellen H. Fukuda},
  year          = {2026},
  eprint        = {2609.30898},
  archiveprefix = {arXiv},
  primaryclass  = {math.OC},
}

@article{merchav2026adaptive,
  author  = {Merchav, Roey and Sabach, Shoham and Teboulle, Marc},
  title   = {The adaptive integrated mean: a novel perspective for acceleration and its extension to {Bregman} first-order methods},
  journal = {Math. Program.},
  year    = {2026},
  doi     = {10.1007/s10107-026-02423-4}
}

@book{lawson1995least,
  author    = {Lawson, Charles L. and Hanson, Richard J.},
  title     = {Solving Least Squares Problems},
  series    = {Classics in Applied Mathematics},
  volume    = {15},
  publisher = {Society for Industrial and Applied Mathematics},
  address   = {Philadelphia},
  year      = {1995},
  doi       = {10.1137/1.9781611971217}
}

@article{slawski2013nonnegative,
  author    = {Martin Slawski and Matthias Hein},
  title     = {Non-negative least squares for high-dimensional linear models: Consistency and sparse recovery without regularization},
  volume    = {7},
  journal   = {Electron. J. Stat.},
  publisher = {Institute of Mathematical Statistics and Bernoulli Society},
  pages     = {3004 -- 3056},
  year      = {2013},
  doi       = {10.1214/13-EJS868}
}

@article{bauschke2003bregman,
  title     = {Bregman monotone optimization algorithms},
  author    = {Bauschke, Heinz H and Borwein, Jonathan M and Combettes, Patrick L},
  journal   = {SIAM J. Control Optim.},
  volume    = {42},
  number    = {2},
  pages     = {596--636},
  year      = {2003},
  publisher = {SIAM}
}

@article{burke1993weak,
  author  = {Burke, James V. and Ferris, Michael C.},
  title   = {Weak Sharp Minima in Mathematical Programming},
  journal = {SIAM J. Control Optim.},
  volume  = {31},
  number  = {5},
  pages   = {1340--1359},
  year    = {1993},
  doi     = {10.1137/0331063}
}

@article{ferris1991finite,
  author  = {Ferris, Michael C.},
  title   = {Finite termination of the proximal point algorithm},
  journal = {Math. Program.},
  volume  = {50},
  pages   = {359--366},
  year    = {1991},
  doi     = {10.1007/BF01594944}
}

@book{powell1981approximation,
  title     = {Approximation theory and methods},
  author    = {Powell, Michael James David},
  volume    = {2148},
  year      = {1981},
  publisher = {Cambridge university press}, 
  address   = {Cambridge}
}

@article{drusvyatskiy2018error,
  title     = {Error bounds, quadratic growth, and linear convergence of proximal methods},
  author    = {Drusvyatskiy, Dmitriy and Lewis, Adrian S},
  journal   = {Math. Oper. Res.},
  volume    = {43},
  number    = {3},
  pages     = {919--948},
  year      = {2018},
  publisher = {INFORMS}
}

@article{fercoq2019adaptive,
  title     = {Adaptive restart of accelerated gradient methods under local quadratic growth condition},
  author    = {Fercoq, Olivier and Qu, Zheng},
  journal   = {IMA J. Numer. Anal.},
  volume    = {39},
  number    = {4},
  pages     = {2069--2095},
  year      = {2019},
  publisher = {Oxford University Press}
}

@article{drusvyatskiy2013tilt,
  title     = {Tilt stability, uniform quadratic growth, and strong metric regularity of the subdifferential},
  author    = {Drusvyatskiy, Dmitriy and Lewis, Adrian S},
  journal   = {SIAM J. Optim.},
  volume    = {23},
  number    = {1},
  pages     = {256--267},
  year      = {2013},
  publisher = {SIAM}
}

\end{document}